\documentclass[11pt,leqno,amscd,amssymb,pstricks,verbatim]{amsproc}
\newtheorem{thm}{Theorem}[section]
\newtheorem{lem}[thm]{Lemma}

\theoremstyle{definition}
\newtheorem{defn}[thm]{Definition}

\newtheorem{prop}[thm]{Proposition}

 \theoremstyle{remark}
\newtheorem{rem}[thm]{Remark}
\newtheorem{cor}[thm]{Corollary}

\numberwithin{equation}{section}

\usepackage{diagrams}

\begin{document}

\def\frakl{{\mathfrak L}}
\def\frakg{{\mathfrak G}}
\def\bbf{{\mathbb F}}
\def\bbl{{\mathbb L}}
\def\bbz{{\mathbb Z}}
\def\bbr{{\mathbb R}}

\def\bvp{\bf{\varphi}}

\def\Der{\mbox{\rm Der}}
\def\Hom{\text{\rm Hom}}
\def\Ker{\text{\rm Ker}}
\def\Lie{\text{\rm Lie}}
\def\id{\mbox {\bf id}}
\def\det{\mbox{\rm det}}
\def\Lie{\mbox {\rm Lie}}
\def\Hom{\mbox {\rm Hom}}
\def\Aut{\mbox {{\rm Aut}}}
\def\Ext{\mbox {{\rm Ext}}}
\def\Coker{\mbox {{\rm Coker}}}
\def\dim{\mbox{{\rm dim}}}

\def\geqs{\geqslant}

\def\ba{{\mathbf a}}
\def\bd{{\mathbf d}}
\def\co{{\mathcal O}}
\def\cn{{\mathcal N}}
\def\cv{{\mathcal V}}
\def\cz{{\mathcal Z}}
\def\cq{{\mathcal Q}}
\def\cf{{\mathcal F}}
\def\cc{{\mathcal C}}
\def\ca{{\mathcal A}}

\def\ggg{{\frak g}}
\def\lll{{\frak l}}
\def\hhh{{\frak h}}
\def\nnn{{\frak n}}
\def\sss{{\frak s}}
\def\bbb{{\frak b}}
\def\ccc{{\frak c}}
\def\ooo{{\mathfrak o}}
\def\ppp{{\mathfrak p}}
\def\uuu{{\mathfrak u}}

\def\p{{[p]}}
\def\modf{\text{{\bf mod}$^F$-}}
\def\modr{\text{{\bf mod}$^r$-}}

\title[Representations of reductive Lie algebras]
{Filtrations in Modular Representations of Reductive Lie Algebras}
\author{Yiyang Li and Bin Shu}
\address{Department of Mathematics, East China Normal University,
Shanghai 200052, China.} \email{52050601014@student.ecnu.edu.cn}
\address{Department of Mathematics, East China Normal University,
Shanghai 200062,  China.} \email{bshu@euler.math.ecnu.edu.cn}
%\thanks{This work is done during the author's visit to Institute of Mathematical Sciences,
%University of Virginia (2004-2005). He would like to express great
%thanks to the department for the hospitality and the financial
%support. This work is supported partially by the NSF of China
%(Grant: 10271047), CSC, a project of science and technology of
%Shanghai, and the Shanghai principal project to Math. Dept. ECNU}

%\subjclass[2000]{14L24, 14L30, 17B10, 17B20, 17B45, 17B50, 20G05,
%20G10, 20G40}

%$\keywords{Nilpotent orbits, support varieties, Frobenius maps,
%${\mathbb F}_q$-rational structures}

\begin{abstract} Let $G$ be a connected reductive algebraic group  over
an algebraically closed field $k$ of prime characteristic $p$, and
$\ggg=\Lie(G)$.
 In this paper, we study representations of the reductive
Lie algebra  $\ggg$ with $p$-character $\chi$ of standard
Levi-form associated with an index subset $I$ of simple roots. With
aid of support variety theory we prove a theorem that  a
$U_\chi(\ggg)$-module is projective if and only if it is a strong
``tilting" module, i.e. admitting  both $\cz_Q$- and
$\cz^{w^I}_Q$-filtrations (to see Theorem \ref{THMFORINV}). Then
by analogy of the arguments in \cite{AK} for $G_1T$-modules, we
construct so-called Andersen-Kaneda filtrations associated  with
each projective $\ggg$-module of $p$-character $\chi$,  and
finally obtain sum formulas from those filtrations.
\end{abstract}

\maketitle

\section{Introduction }
Assume that $k$ is an algebraically closed field of prime
characteristic $p$.
  Let $G$ be a connected reductive algebraic group over $k$,
  and $\ggg = \Lie(G)$.
 Associated with any given linear form $\chi$ on $\ggg$,  $U_{\chi}(\ggg)$ is defined to be the quotient of
   the universal enveloping algebra  $U(\ggg)$ by the ideal generated by all $x^p-x^{[p]}-\chi(x)^p$
    with $x\in \frak{g}$. Each class of irreducible representations of $\ggg$ correspond to a $p$-character
    $\chi$ and the representation theory of $\ggg$ with this $p$-character is the
     ``union" of the representation theory of  $U_{\chi}(\ggg)$. Furthermore, when we restrict
the prime characteristic of the basic fields to the case which we
call ``very good", a well-known result shows that there is a Morita
equivalence between $U_{\chi}(\ggg)$-module category and
$U_{\chi_n}(\frak l)$-module category, where $\frak l$ is a
certain reductive subalgebra of $\ggg$ and $\chi_n$ is a
 so-called nilpotent character of $\lll$ (cf. \cite{KW} and \cite{FP1}).
 This important result enables us to consider the
representations of $U_{\chi}(\ggg)$ just with nilpotent $\chi$. In
the last decade, much progress in modular representations of
reductive Lie algebras have been made. Nevertheless, many  basic
problems remain unsolved.
%We only have a good understand when $\chi$ has standard Levi form (cf. \cite{Jan1})
%and $\chi$ is subregular nilpotent (cf. \cite{Jan4}\cite{Jan5} and
% \cite{Jan6}).%???????????

In this paper, we will focus our concern on the case when  $\chi$
is of standard Levi form which is associated with a subset $I$ of
simple roots of $\ggg$, this means that $\chi$ is regular
nilpotent on the Levi factor $\ggg_I$ for $I$, and is evaluated
$0$ elsewhere.
%(to see \S \ref{SLF} for the definition).
Owing to the work of Friedlander-Parshall and Jantzen
   (cf. \cite{FP1}, \cite{Jan1} and \cite{Jan2}),
    we have a precise classification of simple $U_{\chi}(\ggg)$-modules by "highest weights".
   There are also many good properties in the representations of $U_{\chi}(\ggg)$
   in this case. Especially,  one can study the graded module category by modulo $I$, analogous to
   graded module category in the restricted case, i.e. $\chi=0$
   (the graded structure essentially arises from $G_1T$-module category in representations of algebraic
   groups, where $G_1$ is the kernel of the first
 Frobenius homomorphism and $T$ is the maximal torus of $G$) (cf. \cite[\S 11]{Jan1} and \cite[II. \S9]{Jan3}).
      %In this graded modules
   %category, the representations of $U_{\chi}(\ggg)$ can be studied much batter.

Recall that in the BGG theory for complex representations of
semi-simple Lie algebras, there is a well-known result that each
projective module has a filtration with sub-quotients isomorphic
to Verma modules. In the restricted module category of
$\ggg=\Lie(G)$, Humphreys first proved that each projective module
admits an analogous filtration with sub-quotients isomorphic to
baby Verma modules ${Z}(\lambda)$, called a $Z$-filtration (cf.
\cite{Hum3}). Furthermore, Cline-Parshall-Scott (CPS for short) in
\cite{CPS} proved that a projective $G_1T$-module admits
$Z^w$-filtratios for an arbitrary given $w$ in the Weyl group
 where ${Z}^w(\lambda)$ is a twisted baby Verma module. In the same
time, CPS also proved that a $G_1T$-module is projective if and
only if it is a tilting module, i.e. admiting both $Z$-filtration
and $Z^*$-filtration  (or to say: admitting both $Z$ and
$Z^{w_0}$-filtration for the longest element $w_0$ in the Weyl
group).

In this paper, we prove a strong version of the above result in the
case when $p$-character $\chi$ is of standard Levi form: a $X\slash
\bbz I$-graded module of $U_\chi(\ggg)$ is projective if and only if
it is a strong ``tilting"  module, this is to say, it admits both
$\cz_Q$-filtration and $\cz'_Q$-filtration, where $\cz_Q$ (resp.
$\cz'_Q$) are some $\ppp_I$-induced (resp. $\ppp'_I$-induced)
modules from projective covers of the baby Verma modules of
$U_\chi(\ggg_I)$. Here $\ppp_I$ and $\ppp'_I$ mean respectively
positive and negative parabolic subalgebras associated with $I$ (to
see Theorem \ref{THMFORINV}, where $\cz'_Q$ has another version
$\cz^{w^I}_Q$ by $w^I$-twist which is like $w_0$-twist $Z^{w_0}$
afore-mentioned). The part of ``necessity" in the statement is an
implication of Jantzen's result (to see Proposition
\ref{THMFORPRO}). Here, we complete the other part with aid of
cohomological support variety theory of restricted Lie algebras.

Consequently, a projective $U_\chi(\ggg)$-module $Q$ is a tilting
module\footnote{The inverse of this statement is not true. The
correct inverse statement is just included in Theorem
\ref{THMFORINV}}, i.e. $Q$ admits both $\cz$-filtration and
$^\tau(\cz^*)$-filtration where $\tau$ is an involutative
automorphism of $\ggg$ associated with $I$ (note: $^\tau(\cz^*)$ has
another version $\cz^{w^I}$, to see Lemma \ref{TILTINGDUAL}). So, we
can adopt the filtration introduced by Andersen and Kaneda in
\cite{AK}, and then obtain the sum formulas (to see Theorems
\ref{SumFor1} and \ref{SumFor2}), which is helpful for us to
understand more on simple modules  and their character formulas of
$U_\chi(\ggg)$. On those characters, Lusztig proposed a hope in
\cite{L}, which is still unsolved.

In \cite{AK}, Andersen and Kaneda constructed a filtration
associated with each projective $G_1T$-module $Q$ from vector
spaces $F_{\lambda}(Q)=\Hom_{G_1T}({Z}(\lambda)^{\tau},Q)$, and
then proved that the filtration
 has a sum formula analogous to the Jantzen filtration's,
  (cf. \cite{AK}). We call such  a filtration
   Andersen-Kaneda filtration.

   In the case when $\chi$ is of standard Levi-form for reductive Lie algebras, Andersen-Kaneda
   filtrations do exist and the corresponding sum formulas can be constructed,
   by analogy of A-K's arguments with some mild modifications. For the
   convenience of readers, we complete the arguments.

\vskip0.5cm \textsc{Acknowledgement}:  This work is partially
supported by NSF and PCSIRT of China. The authors express thanks to
the referee for his/her pointing out some false statements in the
original manuscript.

\section{Preliminaries}

%2.1
\subsection{Assumptions}
Throughout this paper, we always assume that $k$ is an
algebraically closed field of prime characteristic $p$.  The
notations will
 generally follow \cite{Jan1}.

 Let $G$ be a connected and reductive algebraic group over $k$, satisfying the following
 three hypotheses
 as in \cite[6.3]{Jan1}:
 \begin{enumerate}
  \item [(H1)] The derived group $\mathcal DG$ of $G$ is simple connected;
  \item [(H2)] The prime $p$ is good for $\ggg$;
  \item [(H3)] There exists a $G$-invariant non-degenerate bilinear form on $\ggg$,
\end{enumerate}
 where $\ggg = \Lie(G)$. Those conditions can be explained as follows: when $\mathcal D G$ is a simple algebraic group,
 then the conditions (H1)-(H3) satisfy if and only if $p$ does not divide
 $n+1$ for type $A_n$; $p>2$ for types $B_n$ $(n\geq 2)$, $C_n$
 $(n>2)$ and $D_n$ $(n\geq 4)$; $p>3$ for type $E_6, E_7$, $F_4$ and $G_2$;
 $p>5$ for $E_8$.

Fix $T$ a maximal torus of $G$. Let $U(\ggg)$ be the universal
enveloping algebra of $\ggg$. For a given  $\chi\in \ggg^*$, set
$U_{\chi}(\ggg)=
 U(\ggg)/J_\chi$, a reduced enveloping  algebra of $\ggg$. Here $J_\chi$ is the ideal of $U(\ggg)$ generated by
$x^{p}-x^{[p]}-\chi {(x)^p}, \mbox{ for all }x\in \ggg $. Set
$X=X(T)$ the character group of $T$, which is a free abelian group
of rank equal to $\dim T$. It contains the subgroup $\mathbb{Z}R$
generated by the root system $R$.

Denote respectively by $R^\pm$ the sets of all positive roots and
all negative roots. For each $\alpha \in R$, let $\ggg_\alpha$
denote the subspace of $\ggg$ corresponding to $\alpha$ and $\frak
n^+=\sum_{\alpha \in R^+}\ggg_\alpha, \frak n^-=\sum_{\alpha \in
R^-}\ggg_\alpha$. We have the triangular decomposition of $\ggg $:
$\ggg=\frak n^+\oplus \frak h\oplus\frak n^-$. Let $\frak
 b^+=\frak h\oplus \frak n^+$ be the Borel subalgebra of $\ggg$ , $\frak h$ is the  Cartan subalgebra
 of $\ggg$.
For each $\alpha\in R$, let $\alpha ^{\vee}$ denote the coroot of
$\alpha$ and $W$ is the Weyl group generated by
 all $s_\alpha$ with $ \alpha \in R $
 and $W_p$ is the affine Weyl group generated by $s_{\alpha,rp} (r\in \mathbb{Z})$ where
  $s_{\alpha,rp}$, $r\in \mathbb{Z}$ are the
affine reflection with
$s_{\alpha,rp}(\mu)=\mu-(\langle\mu,\alpha^{\vee}\rangle-rp)\alpha$.
Define $w.\lambda=w(\lambda+\rho)-\rho,w\in W$, the dot action of
$w$ on $\lambda$ where $\rho$ is half the sum of all positive
roots.

Call $\xi\in\ggg^*$ nilpotent if $\xi$ is in the coadjoint
$G$-orbit of $\chi$ with $\chi(\bbb^+)=0$. By Kac-Weisfeiler and
Friedlander-Parshall's result, up to Morita equivalence, the study
of $U_\xi(\ggg)$-module can be reduced to the study of
$U_{\xi_0}(\ggg_0)$-module for a reductive Lie algebra
$\ggg_0=\Lie(G_0)$ for some connected reductive algebraic group
$G_0$ satisfying Conditions (H1)-(H3), and nilpotent $\xi_0\in
\ggg^*_0$. Thus, we only need study the module category of
$U_\chi(\ggg)$-modules with $\chi$'s being nilpotent, up to Morita
equivalence. In the whole paper, we always assume that
$\chi(\bbb^+)=0$, i.e. $\chi$ is nilpotent.

%2.2

Note that since $\chi(\bbb^+)=0$, any simple $U_0(\hhh)$-module is
one-dimensional $k_\lambda=k$, with $h\cdot 1=\lambda(h)$ for any
$h\in \hhh$ and $\lambda\in \Lambda:=\{\lambda \in \hhh^*\mid
\lambda(h)^p= \lambda(h^{[p]})\}$ which is equal to $X\slash pX$.
And $k_\lambda$ can be extended a $U_0(\bbb^+)$-module with trivial
$\nnn^+$-action. Hence, we have
  an induced module
$${Z}_{\chi}(\lambda)=U_{\chi}(\ggg)\otimes_{U_{0}(\frak b^+)}
k_{\lambda}$$ which is called a baby Verma module. Then each
simple $U_{\chi}(\ggg)$ module is the homomorphic image of some
baby Verma module ${Z}_{\chi}(\lambda),\lambda\in\Lambda$.

\subsection{Standard Levi forms}\label{SLF} We say a
$p$-character $\chi$ has standard Levi form if $\chi$ is nilpotent
and if there exists a subset $I$ of the set of all simple roots
such that

\begin{equation}
 \chi (\ggg _{-\alpha})=\begin{cases}
      \neq 0 & \text{if $\alpha \in{I} $,}\\
      0  & \text{if $\alpha \in{R^+ \backslash I}$}.
    \end{cases}
\end{equation}

As in \cite[\S 10.4; \S 10.5]{Jan1}, when $I$ is the full set of
all simple roots, we call $\chi$ a regular nilpotent element in
$\ggg^*$. When $I=\{0\}$, we know $U_{\chi}(\ggg)=U_{0}(\ggg)$ is
the restricted enveloping algebra of $\ggg$. We denote $R_I$ be
the root system corresponding to simple root set $I$ and
 let $W_I$ is the Weyl group generated by all the $s_\alpha$ with $\alpha \in I$.
Set $w_I$ to be the longest element in $W_I$ and $w_0$ is the
longest element in $W$. We denote
 $w^I=w_Iw_0$.
Denote $\ggg_I=\frak h\oplus \sum_{\alpha\in R_{I}}\frak
 g_{\alpha}$
; $\ppp=\ggg_I+\oplus\uuu^+$; and $\frak p'=\ggg_I+\oplus\uuu^-$,
where $\uuu^+=\sum_{\alpha\in R^+\backslash R^+_I}\ggg_\alpha$,
$\uuu^-=\sum_{\alpha\in R^+\backslash R^+_I}\ggg_{-\alpha}$ are the
nilpotent radicals of $\ppp$ and  of $\ppp'$ respectively.

\vskip0.3cm
 {\sl{From now on, we will always assume that $\chi$ is a
given $p$-character of standard Levi-form, associated with a subset
$I$ of simple roots.}}

\subsection{The category of $X\slash \bbz I$-graded
modules} \label{GRADEPI}

 We are going to study certain $X\slash
\bbz I$-graded $U_\chi(\ggg)$-module category, denoted by $\cc$.
It is defined as follows: Each $U_\chi(\ggg)$-module is a direct
sum of weight spaces of $\hhh$ (note: all weights belong to
$X\slash pX\subset \hhh^*$ as $\chi(\hhh)=0$). If $V\in \cc$, then
each graded component $V_{\lambda+\bbz I}$ with $\lambda\in X$ is
an $\hhh$-submodule, thereby decomposes into weight space for
$\hhh$. Furthermore, a finite-dimensional $X\slash \bbz I$-graded
$U_\chi(\ggg)$-module $V$-belongs to $\cc$ if and only if all
weights of $\hhh$ on $V_{\lambda+\bbz I}$ have the form $d\mu$
with $\mu\in \lambda +\bbz I +pX$ for all $\lambda$. Here
$d\lambda$ means the differential of $\lambda\in X(T)$, which
satisfys $d\lambda\in \Lambda$.

We call a $U_\chi(\ggg)$-module $M$ gradable if there is on $N\in
\cc$ such that $M\cong \cf(N)$. Here $\cf: \cc\rightarrow $
$U_\chi(\ggg)$-module category means the forgetful functor.

\begin{lem} (cf. \cite[1.4]{Jan2})\label{GRADE}\begin{itemize}
\item[(1)] Each simple $U_\chi(\ggg)$-module is gradable.
\item[(2)] Each baby Verma module is gradable. \item[(3)] Each
projective $U_\chi(\ggg)$-module is gradable.\qed
\end{itemize}
\end{lem}

The definition of the $X\slash \bbz I$-graded module category
$\cc$ can be reformulated (and then extended) as follows.

%(cf. \cite[\S 11]{Jan1} and \cite[\S 3]{Jan2})
 Let $U= U(\ggg)/(x_\alpha^p-\chi(x_\alpha)^p, \alpha \in R)$, the PBW theorem ensure
us the following isomorphism:
$$U_{\chi}(\frak n^-)\otimes U(\frak h)\otimes U_{\chi}(\frak n^+)\simeq
U$$
 We shall denote respectively by $U^-, U^0$ and $U^+$ the images of
$U_\chi(\nnn^-), U(\frak h) $ and $ U_\chi(\nnn^+)$ in $U$; these
images is respectively isomorphic to  $U_\chi(\nnn^-),U(\frak h) $
and $ U_\chi(\nnn^+)$. We have a $\mathbb ZR/\mathbb ZI-$grading
on $U$ such that each $\ggg_\alpha$ with $\alpha \in R\cup \{0\}$
is contained in the homogeneous part of degree $\alpha +\mathbb
ZI$. We denote by the homogeneous parts of $U$ by $U_\nu$ with
$\nu\in \mathbb ZR/\mathbb ZI$.

Let $A$ be a noetherian commutative $k$-algebra. Let
$$\pi :U^0\rightarrow A $$
denote a $k$-algebra homomorphism. We call a $U\otimes A$-module $M$
to be $\hhh$-diagonalizable if $$M=\bigoplus_{\phi :\frak
h\rightarrow A}M^\phi,M^\phi:=\{m\in M\mid h.m=\phi(h)m\}.$$ In this
case, we  call $\phi$ an $\frak h$-weight of $M$ if $M^\phi\neq 0$.
Now we can define a category of
 $U\otimes A$-modules $\mathcal C_A$. The objects of $\mathcal C_A$ are
  $U\otimes A$-modules $M$ which are $\hhh$-diagonalizable, together with an $X/\mathbb ZI$-grading
   $M=\bigoplus_{\nu\in X/\mathbb ZI}M_\nu$ of $M$ satisfying the following
   conditions:
\begin{enumerate}
\item[(A)] $M$ is finitely generated over $A$; \item[(B)] $A$
preserves the grading of $M$: $$AM_{\mu+\mathbb{Z}I}\subset
M_{\mu+\mathbb{Z}I};$$ \item[(C)] $U_{\nu+\mathbb ZI}$ shifts the
grading of $M$ by $\nu+\mathbb ZI$:
$$U_{\nu+\mathbb{Z}I}\cdot M_{\mu+\mathbb{Z}I}\subset M_{\mu+\nu+\mathbb{Z}I};$$
\item[(D)] All $\frak h$-weights $\phi$ of $M_{\mu+\mathbb{Z}I}$
have the form:
$$\phi=\pi+d(\mu+\mu'),\mu'\in\mathbb{Z}I.$$
\end{enumerate}

A morphism between objects in $\mathcal C_A$ is a $U\otimes
A$-homomorphism
 preserving the gradings.
 So we can get an induced module in $\mathcal C_A$
 $$\hat{{Z}}_{A}(\lambda)=
 U\otimes_{U^0 U_\chi(\frak n^{+})}A_\lambda$$
 where $A_\lambda$ denotes the $(\frak b^{+})\otimes A$-module $A$ for each
$h\in \frak h $ acting as multiplication by $\pi(h)+d\lambda(h)$.
 %Just
%as the same as in the situation of non-grading,
%$\hat{{Z}}_{A}(\lambda)$
%has a simple head and a simple socle when $\chi$ has standard Levi form.

 We will specialize our choice of the $k$-algebra $A$. When choosing $A=k[t]_{(t)}$
 the localization of the polynomial ring $k[t]$ in one
variable at the maximal ideal generated by $t$, we will denote by
$\tilde {A}$ the fraction field of $k[t]$. The corresponding
category of $U\otimes A$-modules and
  $U\otimes \tilde A$-modules will be denoted by $\cc_A$ and $\cc_{\tilde A}$.
So in the case $A=k[t]_{(t)}$, $\pi:U^0\rightarrow A$ arises from
an algebra homomorphism $\pi_0:U^0\rightarrow k[t]$ which is
defined via sending $h_\alpha$ to $c_\alpha t$ for $c_\alpha \in
k$ with all $\alpha \in R$ satisfying $c_{w_I\alpha}=c_\alpha$ and
$c_\alpha\ne 0$ if and only $\alpha \notin R_I$ (cf.
\cite[3.1/2/9]{Jan2} or \cite[13.2]{Je}).

  When $A=k$, the corresponding category is just $\cc$.

In the sequel,  we will fix an element $\chi$ in $\ggg^{*}$ of
standard Levi form, associated with a subset $I$ of the simple
root system. Then we may define an order relation $\leq$ on
$X/{\mathbb ZI}$
   such that $\mu+\mathbb ZI\leq\lambda+\mathbb ZI$ if and only if there exist integers $m_{\alpha}\geq 0$ with
   $\lambda-\mu=\sum_\alpha m_\alpha\alpha+\mathbb ZI$.
 In category $\mathcal C$, the baby Verma module
 is $$\hat{{Z}}_{\chi}(\lambda)=
 U\otimes_{U^0 U_\chi(\frak n^{+})}k_\lambda,$$
 which satisfies $\mathcal F(\hat{Z_{\chi}}(\lambda))\cong Z_{\chi}(d\lambda)$.
 When $A'=A/tA=k$, we have
 $\hat{{Z}}_{A}(\lambda)\otimes_Ak\simeq\hat{{Z}}_{\chi}(\lambda)$.

%2.5
\subsection{Twists} \label{TWIST} (cf. \cite[\S 3.3]{Jan2}) Let $w$ be an
element of
 $$W^I:=\{w\in W\mid w^{-1}(\alpha)>0,\mbox{ for all }\alpha \in I\}$$
 We can see that $w^I\in W^I$. The PBW basis theorem give us an isomorphism
$$U_{\chi}(w\frak n^-)\otimes U^0\otimes U_{\chi}(w\frak n^+)\simeq
U$$ We can define new categories. In the new category, we get an
induced module for $w\in W^I$
$$\hat{{Z}}_{A}^w(\lambda)=U\otimes_{U^0 U_\chi(w\frak n^{+})}A_\lambda$$
where $A_\lambda$ denotes the $(w\frak b^{+})\otimes A$-module $A$
for each $\frak h $ acting as multiplication by $\pi(\frak
h)+(d\lambda)(\frak h)$.
%Analogous to the non-twist case,
%$\hat{{Z}}_{A}^{w}(\lambda)$ also has a simple head denoted by
%$\hat L^w_\chi(\lambda)$.

When $A=k$, the (twisted) induced module $\hat{{Z}}_{A}^w(\lambda)$
will be denoted by $\hat{{Z}}_\chi^w(\lambda)$. We have the
following facts.

%2.7
\begin{lem} \label{BASICISO} Let $\lambda,\mu\in X$. Then
\begin{enumerate}
\item[(1)] (cf. \cite[Prop. 11.9]{Jan1})
$\hat{{L}}_{\chi}(\lambda)\simeq\hat{{L}}_{\chi}(\mu)\Longleftrightarrow
    \hat{{Z}}_{\chi}(\lambda)\simeq\hat{\mathcal{
     Z}}_{\chi}(\mu)\Longleftrightarrow\mu\in W_{I,p}\cdot\lambda.$
\item[(2)] For $w\in W^I$ there is an equivalence of categories that takes
$\hat{Z}_\chi^w(\lambda)$ to
$\hat{Z}_{w^{-1}\chi}(w^{-1}\lambda)$.
%$$%\hat{{L}}_{\chi}^{w}(\lambda)\simeq\hat{{L}}_{\chi}^{w}
   %(\mu)\Longleftrightarrow
 % \hat{{Z}}_{\chi}^{w}(\lambda)\simeq\hat{Z}_{\chi}^{w}(\mu)\Longleftrightarrow w^{-1}(\lambda)\in
  % W_{wI,p}\cdot w^{-1}(\lambda).$$
\end{enumerate}
\end{lem}

\begin{proof}
(2)  Recall that associated with $w\in W$, there is an automorphism
of $\ggg=\hhh+\sum_{\alpha\in R}\ggg_\alpha$ which stabilizes $\hhh$
and makes $\ggg_\alpha$ into $\ggg_{w\alpha}$. We denote it by $\bar
w$. For $w\in W^I$, define $w^{-1}\chi\in \ggg^*$ via evaluating it
$\chi(\bar wx)$ at $x$ (similarly, we can define $w^{-1}\lambda\in
X$). %Note that $w^{-1}\chi$ is of standard Levi form associated with
%the subset $wI$ of simple roots in the sense of the new simple root
%system twisted by $w$.
Under such an automorphism $\bar w$ for $w\in
W^I$, there is an algebra isomorphism between $U_{w^{-1}\chi}(\ggg)$
and $U_{\chi}(\ggg)$, which gives rise to an category equivalence
between the $X\slash \bbz wI$-graded module category of
$U_{w^{-1}\chi}(\ggg)$ and the $X\slash \bbz I$-graded module
category of $U_\chi(\ggg)$, sending
$\hat{Z}_{w^{-1}\chi}(w^{-1}\lambda)$ to $\hat{Z}_\chi^w(\lambda)$.
\end{proof}

\subsection {$\cz_Q$-$\cz_Q^{w^I}$-Filtrations} \label{QFILTR}
 As in the previous section, we maintain the assumption that
 the $p$-character $\chi$ is a given standard Levi form, in
 connection with a subset of simple roots  $I=\{\alpha\in R
\mid\chi(x_{-\alpha})\neq 0\}$. As the subspaces $\ppp$ and
$\ppp'$ of $\ggg$ are homogeneous, the algebra $U_\chi(\ppp)$ and
$U_\chi(\ppp')$ have a natural grading by $X\slash \bbz I$. When
we extend a $U_\chi(\ggg_I)$-module to a $U_\chi(\ppp)$-module or
a $U_\chi(\ppp')$-module, we will regard  $M$ as a graded module
with $M_0=M$ and $M_\lambda=0$ if all $\lambda \ne 0$. The
situation for $\ppp'$ is the same as $\ppp$.

For each $U_\chi(\ggg_I)$-module $M$ set
$$\cz(M)=U_\chi(\ggg)\otimes_{U_\chi(\ppp)}M$$
and
$$\cz'(M)=U_\chi(\ggg)\otimes_{U_\chi(\ppp')}M.$$

By the arguments in \cite[1.16]{Jan2}, we know that both $\cz(M)$
and $\cz'(M)$ have natural $X\slash \bbz I$-gradings, identifying
$\cz(M)_0$ with $M$ as a $U_\chi(\ggg_I)$-module,
$\cz(M)_\lambda=0$ implies $\lambda\leq 0$ and $\cz'(M)_\lambda\ne
0$ implies $\lambda\geq 0$. The corresponding modules in $\cc$ are
distinguished by wearing  a cap, like $\hat{\cz}(M)$ and
$\hat{\cz'}(M)$.

 Consider ${{Z}}_{\chi,I}(\lambda)=U_{\chi}(\ggg_I)
\otimes_{U_{\chi}(\ggg_I\bigcap\frak
b^+)}k_{\lambda},\lambda\in\Lambda_{\chi} $. Then $\chi|_{\ggg_I}$
is regular nilpotent. By \cite[4.2/3]{FP1}, we know it's a simple
$U_{\chi}(\ggg_I)$ module (also refer to  \cite[\S 10 and \S
11]{Jan1}). Let $Q_{\chi,I}(\lambda)$ be the projective cover of
the $U_{\chi}(\ggg_I)$-module
 ${{Z}}_{\chi,I}(\lambda)$. Thus, we have induced modules of $\ggg$:
$$\cz(Z_{\chi,I}(\lambda))\cong Z_\chi(\lambda).$$
and
$$\cz(Q_{\chi,I}(\lambda)).$$
We denote both by $\cz(\lambda)$ and $\cz(Q,\lambda)$
respectively. By Lemma \ref{GRADE}, there are  corresponding
modules in $\cc$, denoted by $\hat{\cz}(\lambda)$ and
$\hat{\cz}(Q,\lambda)$.

Since
 ${Q}_{\chi,I}(\lambda)$ has a filtration of length $|W_I.\bar\lambda|$
 with all quotient of subsequent terms in the filtration isomorphic to
  ${{Z}}_{\chi,I}(\lambda)$ (cf. \cite[\S10.10]{Jan1}),
  $\cz(Q_{\chi,I}(\lambda))$
 has a filtration of length $|W_I.\bar\lambda|$ with each quotient of
subsequent terms in the filtration isomorphic to
 $Z_{\chi}(\lambda)$. Here and further, $\bar\lambda$ stands for the
 image of $\lambda$ in $X\slash pX$.

Denoting by $ \hat{Q}_{\chi}(\lambda)$  the projective cover of the
simple module $ \hat {L}_{\chi}(\lambda)$, we know  that $
\hat{Q}_{\chi}(\lambda)$ has a filtration with each quotient of
subsequent terms isomorphic
 to $\hat{\cz}(Q,\mu)$ for some $\mu\in X$, while the number of factors isomorphic to a given
 $\hat{\cz}(Q,\mu)$ is equal to $[\hat{{Z}}_{\chi}(\mu):\hat{{L}}_{\chi}(\lambda)]$
 (cf. \cite[Proposition 10.11]{Jan1}). Consequently, the projective module
$\hat{Q}_{\chi}(\lambda)$ has a
 filtration where all factors of subsequent terms are
  isomorphic to $\hat{Z}_\chi(\lambda)$ for some $\lambda\in X$,
  %with
  %all quotient of subsequent terms in the filtration isomorphism to
 %${{Z}}_{\chi}(\mu)$,
 The number (denoted by
 $(\hat{Q}_\chi(\lambda):\hat{Z}_\chi(\mu))$)
  of factors in such a filtration of $\hat{Q}_\chi(\lambda)$ isomorphic to a given $\hat{Z}_{\chi}(\mu)$ is equal to
  $|W_I.\bar\mu|\cdot [\hat{{Z}}_{\chi}(\mu):\hat{{L}}_{\chi}(\lambda)]$.

 Recall that $\ggg_I=\frak h\oplus \sum_{\alpha\in R_{I}}\frak
 g_{\alpha}$; $\ppp=\ggg_I\oplus \uuu^+$. Then
$w^I(\ppp)=\ppp'=\ggg_I\oplus\sum_{\alpha\notin w^I(\bbz
I)}\ggg_{-\alpha}$. Actually, $w^I(I)=-w_I(I)\in \bbz I$, and
$w^I(\alpha)=w_0(\alpha)$ for $\alpha\in R^+\backslash R^+_I$, and
then $w^I(\ggg_I)=\ggg_I$. Hence we can define the twist induced
modules as follows:
$$\cz'(Z_{\chi,I}(\lambda))=U_\chi(\ggg)\otimes_{U_\chi(w^I(\ppp))}Z_{\chi,I}(\lambda)\cong Z^{w^I}_\chi(\lambda).$$
and
$$\cz'(Q_{\chi,I}(\lambda))=U_\chi(\ggg)\otimes_{U_\chi(w^I(\ppp))}Q_{\chi,I}(\lambda).$$
Both of them  will be denoted by $\cz^{w^I}(\lambda)$ and
$\cz^{w^I}(Q,\lambda)$ respectively. The corresponding graded
modules in $\cc$ are denoted by $\hat{\cz}^{w^I}(\lambda)$ and
$\hat{\cz}^{w^I}(Q,\lambda)$.

%  Similarly, let $\hat{L}^{w^I}_{\chi}(\lambda)$ be the simple head of
 % $\hat{Z}^{w^I}_{\chi}(\lambda)$ and  $\hat Q_{\chi}^{w^I}(\lambda)$
  % the projective cover of  $\hat{L}^{w^I}_{\chi}(\lambda)$.
   %By the same argument as above, we know that
    %the projective module $\hat Q_{\chi}^{w^I}(\lambda)$ has a filtration with
  %all subsequent terms in the filtration isomorphic to
  %$\hat{\cz}^{w^I}(Q,\mu)$ and the number of factors isomorphic to a given
   %$\hat{\cz}^{w^I}(Q,\mu)$ is equal to
   %$[\hat{{Z}}^{w^I}_{\chi}(\mu):\hat{{L}}^{w^I}_{\chi}(\lambda)]$.
   %Consequently,  the projective module $\hat Q_{\chi}^{w^I}(\lambda)$ has a filtration with
  %all subsequent terms in the filtration isomorphic to
  %$\hat{Z}_\chi^{w^I}(\mu)$ and the number (denoted by
 %$(\hat{Q}_\chi^{w^I}(\lambda^{w^I}):\hat{Z}^{w^I}_\chi(\mu^{w^I}))$) of factors isomorphic to a given
  % $\hat{Z}_\chi^{w^I}(\mu)$ is equal to
  %$|W_I. \mu|[\hat{{Z}}^{w^I}_{\chi}(\mu):\hat{{L}}^{w^I}_{\chi}(\lambda)]$.

%More generally, for any $w\in W^I$, the $\ggg_I$-module
%$Z_{\chi,I}(\lambda)$ can be regarded a $w(\ggg_I)$-module. In
%this sense, we have $\hat{Z}_\chi^w(\lambda)=$

  \begin{defn}
\begin{itemize}
\item[(1)] Let $M\in \cc$. $M$ is said to has a $\hat{\cz}$-filtration (resp.
$\hat{\cz}^{w^I}$-filtration)
  if there is a filtration with
   each quotient of subsequent terms in the filtration isomorphic
   to $\hat{\cz}_\chi(\mu)$ (resp.
 $\hat{\cz}^{w^I}_\chi(\mu)$ for some $\mu \in X$.
\item[(2)] Let $M\in \cc$. $M$ is said to has a $\hat{\cz}_Q$-filtration (resp.
$\hat{\cz}_{Q}^{w^I}$-filtration)
  if there is a filtration with
   each quotient of subsequent terms in the filtration isomorphic
   to $\hat{\cz}(Q,\mu)$ (resp.
 $\hat{\cz}^{w^I}(Q,\mu)$) for some $\mu \in X$.
\item[(Note:] the same notions in the category $\cc_A$ can be defined,
in the same sense of (1) and (2) respectively.)
 \end{itemize}
\end{defn}

One of Jantzen's results \cite[Proposition 2.1]{Jan2} implies the
following proposition.

\begin{prop}\label{THMFORPRO} Let $\chi\in \ggg^*$ be of standard Levi form.
 If the object $M$ in the category $\mathcal{C}$
is projective, then  $M$ has both $\hat{\cz}_{Q}$-filtration and
$\hat{\cz}_{Q}^{w^I}$-filtration, thereby has both a  $\hat
{\cz}_{\chi}$-filtration and a $\hat{\cz}_{\chi}^{w^I}$-filtration.
\qed
\end{prop}
\begin{proof} Assume $M\in \cc$ is projective, then it's a projective
$U_\chi(\frak p)$-module. By \cite[Proposition 2.1]{Jan2}, $M$ has a
$\hat{\cz}_Q$-filtration. Symmetrically, $M$ has a
$\hat{\cz}^{w^I}_Q$-filtration.
\end{proof}

\vskip0.2cm
\begin{rem}\label{PROJ} \begin{itemize}
\item[(1)]Generally speaking, it's  no longer true that $M$ in $\cc$ is
projective if $M$ admits both $\hat{\cz}_\chi$- and
$\hat{\cz}^{w^I}_\chi$-filtrations. An obvious counter-example is
$M=\hat{Z}_\chi(\lambda)$ for a regular nilpotent $\chi$. In such a
case, $I$ is just  the whole simple roots, and a baby Verma module
coincides with its $w^I$-twist. However, $M$ is projective only when
$\lambda$ is a Steinberg weight (cf. \cite{FP4}).
\item[(2)] By the same argument as \cite[4.19]{AJS}, we have the
following basic facts:
\begin{itemize}
\item[(i)] There is a unique projective module up to isomorphism,
$\hat{Q}_A(\lambda)\in \cc_A$ satisfying
$\hat{Q}_A(\lambda)\otimes_A k\cong \hat{Q}_\chi(\lambda)$, and
$(\hat{Q}_A(\lambda):\hat{Z}_A(\mu))=(\hat{Q}_\chi(\lambda):\hat{Z}_\chi(\mu))$.

\item[(ii)] By the argument in \S\ref{QFILTR}, (i) gives us
 $(\hat{Q}_A(\lambda):\hat{Z}_A(\mu))=|W_I.\bar\mu|\cdot[\hat{Z}_\chi(\mu):\hat{L}_\chi(\lambda)]$.
\item[(iii)] Any projective module $\hat Q$ in $\cc_A$ is isomorphic to a
direct sum of certain $\hat{Q}_A(\lambda)$. Thus, we easily know
that the rank of the free $A$-module
$\Hom_{\cc_A}(\hat{Q},\hat{Z}_A(\mu))$ equal to $(\hat{Q}:
\hat{Z}_A(\mu))$.
\end{itemize}
\end{itemize}
\end{rem}

%2.9
\subsection{Duality} (cf. \cite[\S 11.4; \S 11.5; \S 11.16]{Jan1})
Jantzen constructed a duality $^\tau(-^*)$ on the category
$\mathcal C$.
%There is a good demonstration of Jantzen's duality,
%given by Jessen (cf. \cite[\S11]{Je})).
By \cite[1.14]{Jan2}, there is an automorphism $\tau$ of $G$
satisfying $\tau(T)=T$, with derivative that acts in the following
way:
 $$\tau(x_\alpha)=x_{-w_I\alpha}.$$
 $$\tau(h_\alpha)=h_{-w_I\alpha}.$$
  Note that $\tau ^2=\id$. It has the properties that $\chi\circ\tau^{-1}=-\chi$
  and $\lambda\circ\tau^{-1} =-w_I(\lambda)$ for all $\lambda\in X$.
   In category $\mathcal C$, we have [1, \S11.6]

 %2.10
 \begin{lem} \label{TILTINGDUAL}
 Let $\mu \in X$. $$\hat{{Z}}_{\chi}(\mu)\simeq
{^{\tau}(\hat{{Z}}^{w^I}_{\chi}(\mu^{w^I})^{*})},\; {\text{where
}} \mu^{w^I}=\mu-(p-1)(\rho-w^I\rho).$$ Here $\rho$ is half the
sum of the positive roots of $R$. $\square$
\end{lem}

\begin{rem}\label{REMDUAL} The $\tau$-duality can be extended to the category $\cc_A$.
  Let M be an object of $\mathcal C_A$.
Define $^\tau M$ to be $\Hom_A(M,A)$, as an $A$-module.
  The $U$-action on $ ^\tau(M)$ is defined by:
  $$u\in \ggg,(u.f)(m):=f(-\tau^{-1}(u).m), m\in M .$$
  In general, if $M\in \cc_A$ is a free $A$-module,then
$^\tau(^\tau M)\cong M$. For another $A$-free module $M'\in \cc_A$,
$\Hom_{\cc_A}(M,M')\cong \Hom_{^\tau\cc_A}(^\tau M',^\tau M)$ as
$A$-free modules (compare \cite[\S 1.6]{AK}). Here $^\tau\cc_A$
means an version of the module category corresponding to $-\chi$,
parallel to $\cc_A$.
\end{rem}

In category $\mathcal C_A$, one readily has an analogy of Lemma
\ref{TILTINGDUAL} by a natural way (cf. \cite[\S 13.6]{Je}):

%2.11
\begin{lem} \label{ATILTINGDUAL} Let $\mu\in X$. $$\hat{{Z}}_A(\mu)\simeq
{^{\tau}(\hat{{Z}}^{w^I}_A(\mu^{w^I})^{*})},\; {\text{where }}
\mu^{w^I}=\mu-(p-1)(\rho-w^I\rho).$$ Here $\rho$ is half the sum of
the positive roots of $R$. $\square$
\end{lem}

\section{Baby Verma modules and their twists}

In this section, we will give some computation on hom-spaces and
extensions between (baby) Verma modules and their twists in the
module category $\mathcal C_A$, which will be used later. Before
that, we first prove  some general formulation of Lemma
\ref{BASICISO} which will be the start point of our argument in
the sequel.

\begin{lem}\label{LEMISO3.8}  Let $\lambda,\mu\in X$. Then
%\begin{enumerate}
%\item[(1)] $\hat{{Z}}_A(\lambda)\simeq\hat{{Z}}_A(\mu)
%\Longleftrightarrow\hat{{L}}_A(\lambda)
%\simeq\hat{{L}}_A(\mu)\Longleftrightarrow\mu\in
%W_{I,p}\cdot\lambda .$ \item[(2)]
%$\hat{{Z}}^{w^I}_A(\lambda)\simeq\hat{{Z}}^{w^I}_A(\mu)
%\Longleftrightarrow\hat{{L}}^{w^I}_A(\lambda)
%\simeq\hat{{L}}^{w^I}_A(\mu)\Longleftrightarrow\mu\in
%W_{I,p}\cdot\lambda .$
%\item[(1)]
$\Hom_{\mathcal{C}_A}(\hat{{Z}}_A(\lambda),\hat{{Z}}_A(\mu))\simeq
A\Longleftrightarrow \mu\in W_{I,p}\cdot\lambda.$
%\item[(2)] More generally, for $w\in W^I$
%$$\Hom_{\mathcal{C}_A}(\hat{{Z}}^{w}_A(\lambda^w),\hat{{Z}}^{w}_A(\mu^w))\simeq
%A\Longleftrightarrow w^{-1}\mu^w\in W_{wI,p}\cdot w^{-1}\lambda^w.$$
%\end{enumerate}
\end{lem}

\begin{proof} First, we assert that all
  $\hat{Z}_{\tilde A}(\lambda)=\hat{{Z}}_A(\lambda)\otimes_{A}\tilde{A}$ are simple in $\cc_{\tilde
  A}$. We will prove this by standard argument. Consider $\ggg_{\bar K}=\ggg_K\otimes_K
  \bar K$, where $K=\tilde A$, $\ggg_K$ the extension of Lie algebra $\ggg$
  by field extension of $K\slash k$ and $\bar K$ is the
  algebraically closure of $K$. Then we have
  \begin{equation}\label{KMOD}\hat{Z}_{\bar K}(\lambda)\cong
  \hat{Z}_{\chi_\pi}(\lambda)
  \end{equation}
   where $\chi_\pi\in \ggg^*_{\bar K}$
  is defined via $\chi_\pi(h_\alpha)=c_\alpha t-(c_\alpha
  t)^{1\over p}\ne 0$ if $\alpha\in R\backslash R_I$,
  $\chi_\pi(h_\alpha)=0$ if $\alpha\in R_I$ and
  $\chi_\pi|_{\nnn^{\pm}}=\chi|_{\nnn^{\pm}}$. Here $c_\alpha$ is
  defined as in \S \ref{GRADEPI}. Then $\chi_\pi$ has a
  Chevalley-Jordan decomposition $\chi_{\pi,s}+\chi_{\pi,n}$ where
    $\chi_{\pi,s}\in \ggg^*_{\bar K}$ is the trivial extension of
  $\chi_{\pi}|_{\hhh}\in \hhh^*\subset \ggg^*_{\bar K}$, and $\chi_{\pi, n}\in \ggg^*_{\bar K}$ is
  the trivial extension of $\chi\in \ggg^*\subset \ggg^*_{\bar
  K}$. Thus the centralizer $\ccc_{\ggg_{\bar K}}(\chi_\pi)$ of $\chi_\pi$ in $\ggg_{\bar K}$
   coincides with $\ggg_I\otimes \bar K$. By
  \cite[2.4]{KW} (more precisely \cite[3.2 and 8.5]{FP1}),   $U_{\chi_\pi}(\ggg_{\bar K})$ is Morita
  equivalent to $U_\chi(\ggg_I\otimes \bar K)$, the latter of which is
  of standard Levi form. All Verma modules of $U_\chi(\ggg_I\otimes \bar
  K)$ are simple, which implies the corresponding baby Verma
  modules $U_{\chi_\pi}(\lambda)$ of $U_{\chi_\pi}(\ggg_{\bar K})$ are simple, under the Morita
  equivalence. So $\hat{Z}_{\bar K}(\lambda)$ is simple (notice (\ref{KMOD})).
  From $\hat{Z}_{\bar K}(\lambda)\cong \hat{Z}_{\tilde A}(\lambda)\otimes_{\tilde A} \bar K$, it follows that
  all $\hat{Z}_{\tilde A}(\lambda)$ are simple.

  Furthermore,  by  the Morita equivalence and Lemma
  \ref{BASICISO} (over $\bar K$), we know
  $$\Hom_{\cc_{\tilde A}}(\hat{Z}_{\tilde A}(\lambda),\hat{Z}_{\tilde A}(\mu))
  =\begin{cases} \tilde A &{\text { if }}
  \lambda \in W_{I,p}\cdot \mu, \cr
  0 &{\text{ otherwise. }} \end{cases}$$
Observe that $\tilde A$ is $A$-flat, and the $A$-free module
$\Hom_{\cc_A}(\hat{Z}_A(\lambda),\hat{Z}_A(\mu))$ satisfies by
\cite[2.38]{CR}
$$\Hom_{\cc_A}(\hat{Z}_A(\lambda),\hat{Z}_A(\mu))\otimes_A \tilde A\cong
\Hom_{\cc_{\tilde A}}(\hat{Z}_{\tilde A}(\lambda),\hat{Z}_{\tilde
A}(\mu)).$$
  Hence   $\Hom_{\cc_A}(\hat{Z}_A(\lambda),\hat{Z}_A(\mu))=A$ if and only if $\lambda\in
  W_{I,p}\cdot \mu$.
\end{proof}

%3.10
%\begin{rem} With Lemmas \ref{LEMISO3.8} and \ref{LEMISO3.9},
%it is easily seen that the argument for the first part (before
%Theorem \ref{THMFORPRO}) in this section is still valid in
%$\mathcal C_A$.
%\end{rem}

The following lemmas are crucially useful in the last section.
%3.11
\begin{lem}\label{LEMISO3.11}
\begin{itemize}
\item[(1)] Let $\lambda\in X$, and
$\lambda^\sigma=\lambda+(p-1)(\sigma\rho-\rho)$ for $\sigma\in W^I$.
Then for $w,w'\in W^I$,
$\Hom_{\cc_A}(Z_A^w(\lambda^w),Z_A^{w'}(\lambda^{w'}))\cong A$.
\item[(2)] In particular, for $\lambda\in X$, and $w^I=w_Iw_0$,
$$\Hom_{\cc_A}(\hat{{Z}}_A(\lambda),
\hat{{Z}}^{w^I}_A(\lambda^{w^I}))\simeq A$$ and
$$\Hom_{\cc_A}(\hat{{Z}}^{w^I}_A(\lambda^{w^I}),\hat{{Z}}_A(\lambda))\simeq A.$$
\end{itemize}
\end{lem}
\begin{proof} We can prove the lemma, following the argument in the proof of \cite[4.7]{AJS} (or to see
\cite[1.7]{AK}).
  \end{proof}

We have the following reformulation of  Lemma \ref{LEMISO3.11},
dealing with general situations.

\begin{lem}\label{HOMZZ}
Let $\lambda,\mu\in X$.  For any $w\in W^I$,
$$\Hom_{\mathcal{C}_A}(\hat{{Z}}_A^w(\lambda^w),\hat{Z}_A(\mu))\simeq
    \left\{
    \begin{array}{ll}
      A & \mbox{if
            $\mu\in W_{I,p}\cdot \lambda$}\\
      0  & \mbox{others}
    \end{array}
    \right.$$
 \end{lem}

\begin{proof} As $\tilde A$ is $A$-flat, applying \cite[2.38]{CR} we know that
$$\Hom_{\mathcal{C}_A}(\hat{{Z}}_A^w(\lambda^w),\hat{Z}_A(\mu))
\otimes_A \tilde A\cong  \Hom_{\cc_{\tilde A}}(\hat{Z}_{\tilde
A}^w(\lambda^w),\hat{Z}_{\tilde A}(\mu)).$$ If
$\Hom_{\mathcal{C}_A}(\hat{{Z}}_A^w(\lambda^w),\hat{Z}_A(\mu))$ is
nonzero, then $\Hom_{\cc_{\tilde A}}(\hat{Z}_{\tilde
A}^w(\lambda^w),\hat{Z}_{\tilde A}(\mu))$ is nonzero. By the
argument in the proof of Lemma \ref{LEMISO3.8}, we know that both
$\hat{Z}_{\tilde A}^w(\lambda^w)$ and $\hat{Z}_{\tilde A}(\mu)$ are
simple. So $\hat{{Z}}_{\tilde A}^w(\lambda^w)\cong \hat{{Z}}_{\tilde
A}(\mu)$. Lemma \ref{LEMISO3.11}(1) tells us that $\hat{{Z}}_{\tilde
A}^w(\lambda^w)\cong \hat{{Z}}_{\tilde A}(\lambda)$. Using Lemma
\ref{BASICISO}, we finally obtain $\mu \in W_{I,p}\cdot \lambda$.
Conversely, if $\mu\in W_{I,p}\cdot \lambda$, by the above argument
we know
$\Hom_{\cc_A}(\hat{Z}_A^w(\lambda^w),\hat{Z}_A(\mu))\otimes_A\tilde
A\cong \tilde A$. Note that the first term in the tensor product is
a free-A module of finite rank, thereby of rank one. We complete the
proof.
\end{proof}

 %Actually,
%Jantzen can concretely construct a generator in
%$\hbox{Hom}_{\mathcal{C_A}}(\hat{{Z}}_A(\lambda),
%\hat{{Z}}^{w^I}_A(\lambda^{w^I}))$ (cf. \cite[\S 3.10]{Jan2}).

Choose a reduced expression  $w^I=s_1s_2\cdots s_N$ where
$s_i=s_{\alpha_i}$ for some simple roots $\alpha_i$. Set
$w_i=s_1\cdots s_{i-1}, i=1,\cdots,N+1$ with convention $w_1=1$.
Then all $w_i\alpha_i$ are distinct, they are exactly the positive
roots made negative by $w_0w_I=(w^I)^{-1}$, constituting
$R^+\backslash R^+_I$. Furthermore, the positive roots made negative
by $w_i^{-1}$ are exactly $w_1\alpha_1,\cdots,w_{i-1}\alpha_{i-1}$.
This shows that $w_i^{-1}(I)\subset R^{+}$, thereby $w_i\in W^I$.
Hence we can construct the $w_i$-twist baby Verma modules
$Z^{w_i}_A(-)$ and $Z^{w_i}_\chi(-)$. Set
$$\hat{{Z}}^i_A={Z}^{w_i}_A(\lambda^{w_i}) \;\;\;\;\;
(\mbox{ resp. }\hat{{Z}}^i_{\chi}={Z}^{w_i}_{\chi}
(\lambda^{w_i})).$$ Then $\hat{{Z}}_A^1=\hat{{Z}}_A(\lambda)$ and
$\hat{{Z}}_A^{N+1}={Z}_A^{w^I}(\lambda^{w^I})$. By  Lemma
\ref{LEMISO3.11}(1), we can take a generator $\varphi_i\in
\Hom_{\cc_A}(\hat{Z}^i,\hat{Z}^{i+1})\cong A$, which is unique up
to units in $A$. Similarly, we take a generator $\varphi'_i\in
\Hom_{\cc_A}(\hat{Z}_A^{i+1},\hat{Z}_A^{i})\cong A$. Precisely,
those module homomorphisms $\varphi_i$ and $\varphi'_i$ can be
taken respectively via: sending $1\otimes 1$ to
$x^{p-1}_{w_i\alpha_i}\otimes 1$, and sending $1\otimes 1$ to
$x^{p-1}_{-w_i\alpha_i}\otimes 1$.

By change of rings $-\otimes_A k$ (note that $k\cong A\slash tA$),
one has extensions of both $\varphi_i$ and $\varphi'_i$:
$\bar\varphi_i=\varphi_i\otimes_A k$ and
$\bar\varphi'_i=\varphi'_i\otimes_A k$ in
$\Hom_\cc(Z^i_\chi,Z^{i+1}_\chi)$ and
$\Hom_\cc(Z^{i+1}_\chi,Z^i_\chi)$ respectively. It can be known from
the forthcoming lemma that if $\langle \lambda+\rho,
w_i\alpha_i^\vee\rangle\equiv 0 \mod (p)$, then both $\bar\varphi_i$
and $\bar\varphi'_i$ are isomorphisms in $\cc$ .

\begin{lem}\label{AK1.9} Up to units in $A$ the following
statements hold
\begin{itemize}
\item[(1)] $\varphi_i\circ \varphi'_i=\begin{cases}
t\id_{Z^{i}_A}, \;\; &{\text{\rm{ if }}}\langle \lambda+\rho,
w_i\alpha_i^\vee\rangle\neq 0 \mod (p), \cr \id_{Z^{i}_A},
&{\text{\rm{ otherwise }}}. \end{cases}$

\item[(2)] $\varphi'_i\circ \varphi_i=\begin{cases}
t\id_{Z^{i+1}_A}, \;\; &{\text{ if }}\langle \lambda+\rho,
w_i\alpha_i^\vee\rangle\neq 0 \mod (p), \cr \id_{Z^{i+1}_A},
&{\text{\rm {otherwise} }}.
\end{cases}$
\item[(3)] If $\langle \lambda+\rho, w_i\alpha_i^\vee\rangle\equiv
0 \mod (p)$, then both $\varphi_i$ and $\varphi'_i$ are
isomorphisms.
\end{itemize}
\end{lem}

\begin{proof} (1) Thanks to Lemma \ref{LEMISO3.11}, it is sufficient
to evaluate $\varphi'_i\circ \varphi_i$ on a single non-zero element
$v$ in $Z_A^{i}$, say $v_0=1\otimes 1$. Recall that for $0\leq s\leq
p-1$, $v_s:={x^s_{-w_i\alpha_i}\over s!}\otimes 1\in Z^{i}_A$ and
$v'_s:={x^s_{w_i\alpha_i}\over s!}\otimes 1\in Z^{i+1}_A$, we have
$$x_{w_i\alpha_i}v_s=(\pi(h_{w_i\alpha_i})+\langle \lambda
+\rho,w_i\alpha_i^\vee\rangle -s)v_{s-1},$$
 where $\alpha_i^\vee$
is the coroot. Note that $\varphi_i(v_0)=(p-1)!v'_{p-1}$ and
$\varphi'_i(v'_0)=(p-1)!v_{p-1}$. Hence we have
\begin{eqnarray}\label{FORDUALB} \varphi'_i\circ
\varphi_i(v_0)&=(p-1)!\prod^{p-1}_{j=1}(\pi(h_{w_i\alpha_i})+\langle
\lambda+\rho, w_i\alpha_i^\vee\rangle+j)v_0\cr
&=(p-1)!\prod^{p-1}_{j=1}(c_{w_i\alpha_i}t+\langle \lambda+\rho,
w_i\alpha_i^\vee\rangle+j)v_0.
\end{eqnarray}
 When $\langle \lambda+\rho,
w_i\alpha_i^\vee\rangle\equiv 0 \mod(p)$, the coefficient above is
a unit in $A$. When $\langle \lambda+\rho,
w_i\alpha_i^\vee\rangle\ne 0\mod(p)$, the coefficient above is
equal to $tu$ for  a unit $u$ in $A$. Hence, we get the first
statement.

The proof of (2) is similar. The third statment is an immediate
consequence of (1) and (2).
\end{proof}

\begin{rem}\label{FORDUALB2} We have general formulas in connection with
(\ref{FORDUALB}) which will be used:
$$\varphi_i(v_s)=(-1)^s {(p-1)!\over
s!}\prod_{j=1}^s(\pi(h_{w_i\alpha_i})+\langle
\lambda+\rho,w_i\alpha_i^\vee\rangle-j)v'_{p-1-s}$$ and
$$\varphi'_i(v'_s)= {(p-1)!\over
s!}\prod_{j=1}^s(\pi(h_{w_i\alpha_i})+\langle
\lambda+\rho,w_i\alpha_i^\vee\rangle+j)v_{p-1-s}$$
\end{rem}

%By \cite[5.13]{AJS}, we have
%\begin{eqnarray} \varphi_i
%{\text{ is bijective on the }}
% \lambda
% {\text{-weight space if }}
%s_1\cdots s_{i-1}\alpha_i\in R^+
%\end{eqnarray}

 We can take the
generator $\varpi$ of
$\hbox{Hom}_{\mathcal{C}_A}(\hat{{Z}}_A^1,\hat{Z}_A^{N+1})$
 as the composite
$$\hat{{Z}}_A^1\stackrel{\varphi_1}{\rightarrow}
\hat{{Z}}_A^2\stackrel{\varphi_2}{\rightarrow}
\cdots\stackrel{\varphi_{N+1}}{\rightarrow}\hat{Z}_A^{N+1}.$$
Similarly we can take the generator
$\varpi'=\varphi'_1\circ\cdots\varphi'_N$ in
$\Hom_{\cc_A}(\hat{Z}_A^{N+1},\hat{Z}_A^1)$. Then we have the
following direct corollary to the above lemma.

\begin{cor} \label{KEYCOR} Keep the assumption and notations. Then
$$\varpi'\circ \varpi=t^{N(I,\lambda)}\id_{\hat{Z}^i_A}\;\;\;\;
{\text{ and }}\varpi\circ
\varpi'=t^{N(I,\lambda)}\id_{\hat{Z}^{i+1}_A},$$ where
$N(I,\lambda)=\#\{\alpha\in R^+\backslash R^+_I\mid \langle
\lambda+\rho, \alpha^\vee\rangle\ne 0 \mod (p)\}$. \qed
\end{cor}

%3.14
\begin{lem} \label{AJSORDER}
 \begin{enumerate}
 \item[(1)]  Let $\lambda, \mu \in X$. If $\Ext^1_{\mathcal{C}_A}(\hat{{Z}}_A(\lambda),
\hat{{Z}}_A(\mu))\neq 0$, then $\mu+\mathbb{Z}I\geq
\lambda+\mathbb{Z}I$. \item[(2)] Assume that a module $M$ in
category $\mathcal{C}_A$ has a $\hat{\cz}$-filtration. Then we can
find a $\hat{\cz}$-filtration of  $M$:
$$M=M_0\supset M_1\supset M_2\supset\cdots \supset M_r\supset M_{r+1}=0$$
with $M_i\slash M_{i+1}\cong \hat{Z}_A(\lambda_i)$, $\lambda_i\in
X$, $i=0,1,\cdots,r$ satisfying that
 $\lambda_i+\bbz I\geq\lambda_j+\mathbb{Z}I$ implies  $i\leq
j.$
\end{enumerate}
\end{lem}

\begin{proof}  In analogy of \cite[Lemma 2.14]{AJS}, we consider the exact sequence  in the category
$\mathcal{C}_A$:
$$0\rightarrow \hat{{Z}}_A(\mu)\rightarrow N
\rightarrow \hat{{Z}}_A(\lambda)\rightarrow0$$ Let $v\in
N_{\lambda}$ be an inverse image of the standard generator
$v_0=1\otimes1$ of $\hat{{Z}}_A(\lambda)$. If $x_{\alpha}v=0$ for
all $\alpha \in R^{+}$, then there is a homomorphism
$\hat{{Z}}_A(\lambda)\rightarrow
 N$ maps $v_0$ to $v$, splitting the above exact sequence. So, if the
 above exact sequence does not split, there must be $\alpha > 0$ with $x_\alpha v\neq0$
 . We get that $\lambda+\alpha$ is a weight of $M$, thereby is a weight of $\hat{Z}_A(\mu)$
 because of the exact sequence. The first statement is proved. The
 second is a consequence of (1), by induction on the length of the
 filtration.
 \end{proof}

%3.15
\begin{rem}\label{rem3.15}
\label{ORDERFORFIL}In Lemma \ref{AJSORDER}(1), if we replace $\hat{{Z}}_A(\mu)$ and $\hat{{Z}}_A(\lambda)$ with
 $\hat{{Z}}^{w^I}_A(\mu^{w^I})$ and $
\hat{{Z}}^{w^I}_A(\lambda^{w^I})$ respectively, we can get the
similar result. In (2), if we replace
$\hat{{Z}}$-filtration with
 $\hat{\cz}^{w^I}$-filtration, we can find that the
 $\hat{\cz}^{w^I}$-filtration of $M$ has the properties:
 if
$\lambda_i+\mathbb{Z}I\geq\lambda_j+\mathbb{Z}I,\mbox{ then }i\geq
j.$

\end{rem}

\section{Projective modues and $\hat{\cz}_{Q}$-$\hat{\cz}_{Q}^{w^I}$-filtrations}

 %3.3

The following theorem shows that the inverse of the statement in
Proposition \ref{THMFORPRO} concerning $\cz_{Q}$-filtrations is
true.

\begin{thm} \label{THMFORINV}  Maintain the assumption as in
Proposition \ref{THMFORPRO}. Then $M$ is projective in the
category $\cc$ if and only if $M$ has both $\hat \cz_{Q}$- and
$\hat\cz_{Q}^{w^I}$-filtrations.
\end{thm}

In order to prove the theorem, we need some knowledge about
support varieties and rank varieties (one can refer to the
definitions \cite{FP1} \cite{FP2} \cite{FP3} and \cite{Jan1}).

%3.4.
 \subsection{Support varieties and rank varieties.}
Given a vector space $V$ over $k$ let $V^{(-1)}$ denote the
$k$-space with underlying
 abelian group $V$ and with scalars $\lambda\in k$ acting through multiplication
  by $\lambda^p$. Let $||\ggg||$ denote the affine algebraic variety associated
  to the commutative Noetherian $k$-algebra $\bigoplus_{i\geq0}\hbox{Ext}^{2i}_{U_0(\ggg)}(k,k)^{(-1)}.$

\begin{defn} (\cite[6.1]{FP1})
 Let $M$ be a finite dimension $U_\chi(\ggg)$-module. The support
 variety $|| \ggg||_M$ of $M$ is the closed subvariety of $||\ggg||$
  given as the support of the $\hbox{Ext}^{*}_{U_0(\ggg)}(k,k)^{(-1)}$-module
  $\hbox{Ext}^{*}_{U_\chi(\ggg)}(M,M)^{(-1)}$.
  \end{defn}

   There is a natural finite morphism $\Phi: ||\ggg||\rightarrow\ggg$
   of affine varieties (cf. \cite[\S 2.1]{FP2}), defined by using the $k$-algebra homomorphism
    $$S(\ggg^{*})\rightarrow\bigoplus_{i\geq0}\hbox{Ext}^{2i}_{U_0(\ggg)}(k,k)^{(-1)}$$
     induced from the natural (Hochschild) map $\ggg^{*}\rightarrow
     \hbox{Ext}^{2}_{U_0(\ggg)}(k,k)^{(-1)}.$ By Jantzen's theorem (cf. \cite[Satz 2.14]{Jan6}),
     $\Phi(||\ggg||)$ identifies with the closed subvariety
     $\mathcal N_p(\ggg)$ of $\ggg$ where $\mathcal N_p(\ggg)=\{x\in \ggg\mid x^{[p]}=0\}$.
          From the morphism $\Phi: ||\ggg||\rightarrow\mathcal N_p(\ggg)$,
     one can get  a Zariski closed conical subset $\Phi(||\ggg||_M)$ in
      $\cn_p(\ggg)$, which can be identified with $\{0\}\bigcup\{x\in \cn_p(\ggg) \mid
      M|_{kx}$ is not a free $kx$-module$\}$ (cf.
      \cite{FP1} and \cite{FP3}). % $\Phi(||\ggg||_M)$ is called the rank variety of $M$.

 %3.5
 \begin{lem} (cf. \cite[Prop 6.2]{FP1}) \label{FPTHM}M is a projective $U_\chi(\ggg)$-module if and
  only if $\Phi(||\ggg||_M)=0.$ $\square$
  \end{lem}

 %3.6
 %\begin{lem}  (Premet's theorem, cf. \cite[Thm 1.1]{Pre}) \label{PRETHM} Assume  $p$ is good
 % for $\ggg$. Let $M$ be a nonzero $U_\chi(\ggg)$-module. Then
  % $\Phi(||\ggg||_M)\subseteq \frak c_{\ggg}(\chi)\cap\mathcal N_p(\ggg)$,
   % where $\frak c_{\ggg}(\chi)=\{x\in \ggg\mid \chi([x,\ggg])=0\}$. $\square$

  %\end{lem}

  Now let us prove Theorem \ref{THMFORINV}.
    \begin{proof} $(\Longrightarrow)$ It is what Proposition
    \ref{THMFORPRO} says.
    %an implication of
    %\cite[Proposition 1.21]{Jan2}.
    %where the argument in \cite[Proposition II.11.2]{Jan3} is needed.

$(\Longleftarrow)$ Suppose $M\in \cc$ has a $\hat
\cz_{Q}$-filtration and $\hat{\cz}_{Q}^{w^I}$-filtration. Observe
that $M$ is projective in $\cc$ if only if $\cf(M)$ is a
projective $U_\chi(\ggg)$-module. By Lemma \ref{FPTHM}, we only
need to prove that $\Phi(||\ggg||_{\cf(M)})=0$. Owing to
\cite[Prop 7.1]{FP1}, we know $\Phi(||\ggg||_{\cf(M)\otimes
\cf(M)})=
 \Phi(||\ggg||_{\cf(M)})\cap \Phi(||\ggg||_{\cf(M)})$. Hence we only need consider
  $\Phi(||\ggg||_{\cf(M)\otimes \cf(M)})$.

  As $M$
  has a $\hat {\cz}_{Q}$-filtration and a $\hat{\cz}_{Q}^{w^I}$-filtration,
   $M\otimes M$ has a filtration where each sub-quotient in the filtration is isomorphic
   to $\hat{\cz}(Q,\lambda_i)\otimes M$ which admits another filtration with quotients of sub-quotients
   isomorphic to $\hat{\cz}(Q,\lambda_i)\otimes\hat{\cz}^{w^I}(Q,\lambda_j)$  for
 some $\lambda_i, \lambda_j \in X$.
 Hence, $\Phi(||\ggg||_{\cf(M)\otimes \cf(M)})\subset \bigcup_{i,j}
 \Phi(||\ggg||_{\cz(Q,\lambda_i)\otimes \cz^{w^I}(Q,\lambda_j)})$ (cf. \cite[page 1085]{FP1}).
Associated with each component in the union, we
have$$\Phi(||\ggg||_{\cz(Q,\lambda_i)\otimes\cz^{w^I}(Q,\lambda_j)})=
 \Phi(||\ggg||_{\cz(Q,\lambda_i)})\cap\Phi(||\ggg
 ||_{\cz^{w^I}(Q,\lambda_j)}).$$

By the same arguments in \cite[Remark 7.5]{FP1}, we have
$\Phi(||\ggg
 ||_{\cz(Q,\lambda_j)})\subset \ppp$ and
$\Phi(||\ggg
 ||_{\cz^{w^I}(Q,\lambda_j)})\subset w^I(\ppp)=\ppp'$.
  Thus,
$$\Phi(||\ggg||_{\cz(Q,\lambda_i)\otimes\cz^{w^I}(Q,\lambda_j)})=
 \Phi(||\ggg||_{\cz(Q,\lambda_i)})\cap\Phi(||\ggg
 ||_{\cz^{w^I}(Q,\lambda_j)})
\subseteq \ppp\cap \ppp'=\ggg_I.$$

On the other hand, both $\cz(Q,\lambda_i)$ and
$\cz^{w^I}(Q,\lambda_j)$ are projective as
$U_\chi(\ggg_I)$-module, thereby both
$\Phi(||\ggg||_{\cz(Q,\lambda_i)})$ and $\Phi(||\ggg
 ||_{\cz^{w^I}(Q,\lambda_j)})$ intersect $\ggg$ at
 $0$, owing to Lemma \ref{FPTHM}.
Hence, $ \Phi(||\ggg||_{\cz(Q,\lambda_i)})\cap\Phi(||\ggg
 ||_{\cz^{w^I}(Q,\lambda_j)})=0$. Thus, we have proved that $\Phi(||\ggg||)_M=0$. Hence, $M$ is
 projective. The proof is completed.
 \end{proof}

%By \cite[Remark 7.5]{FP1} we know that
% $\Phi(||\ggg||_{\hat{\mathcal Z}_\chi(\lambda_i)})\subset \frak n^+$.

%By the arguments in \cite[11.17 and 11.18]{Jan1}),  we know  that
%$$\hat\cz_{Q,\chi}^{w^I}(\lambda_j)\cong
%{^\tau\hat\cz_{Q,\chi}_{-\chi}(-(w^I)^{-1}\lambda_j)}.$$
 % In the meantime,
%it's readily proved that $\Phi(||{^\tau M}||)={\tau\Phi(||M||)}$.

%Combining Premet's Theorem (Lemma \ref{PRETHM}), we have
% $$\Phi(||\ggg||_{\hat{\mathcal Z}_\chi(\lambda_i)\otimes\hat{\mathcal Z}_{\chi}^{w^I}(\lambda_j)})=
% \Phi(||\ggg||_{\hat{\mathcal Z}_\chi(\lambda_i)})\cap\Phi(||\ggg
% ||_{\hat{\mathcal Z}^{w^I}_{\chi}(\lambda_j)})
 %$\Phi(||\ggg||_{\hat{\mathcal Z}_\chi(\lambda_i)})\cap\Phi(||\ggg
 %||_{\hat{\mathcal Z}^{w^I}_{\chi}(\lambda_j)})
% \subseteq c_{\ggg}(\chi)\cap \ggg_I.$$

%\begin{rem} \begin{enumerate}
%\item[(1)] When $\chi$ is the regular nilpotent, this theorem is
%not true in general. It is because at this case $w^I$ is $1$ and
%$M$ has only $\cz_{\chi,I}$-filtration. If $M$ has just only
%$\cz_{\chi,I}$-filtration, we know it is not projective in general
%(except that $\lambda$ is special (cf. \cite[\S 4.12]{Jan2}).
%\item[(2)] When $\chi$ is zero, this theorem is well-known as in
%\cite[II \S 11.4]{Jan3} which is first proved in \cite{CPS} and in
%\cite[prop 2.5]{AK} Andersen and Kaneda give it a another proof.
%So we here provide a new proof even in this special case.
%\end{enumerate}
%end{rem}

%3.16
We immediately have the following corollary.

\begin{cor}\label{PROJC} If $M$ is a projective module in $\cc$, then $^\tau M$
must be projective in $^\tau\cc$, where $^\tau\cc$ is the $X\slash
\bbz I$-graded module category of $U_{-\chi}(\ggg)$. In particular,
in this case $^\tau M$ admits a filtration with filtration quotient
factors $^\tau\hat{Z}(\lambda)$ for some $\lambda \in X$.\qed
\end{cor}

\begin{rem}\label{PROJCORA} By Remarks \ref{PROJ}(2) and \ref{REMDUAL}£¬ the $\cc_A$-version of
Corollary \ref{PROJC} is also true. This is to say, if $\hat Q$ is
projective in $\cc_A$, then $^\tau \hat Q$ admits a filtration with
filtration quotient factors like $^\tau\hat{Z}_A(\lambda)$.
\end{rem}

\begin{prop}(Compare with \cite[Proposition 4.11]{AJS})
\label{AJS411} Let $\chi$ be of standard Levi form associated with a
proper subset $I$ of simple roots. Then
$$\Ext^n_{\cc}(\hat{\cz}(Q,\lambda),
\hat{\cz}^{w^I}(Q,\mu^{w^I}))=0$$ for all $\lambda,\mu \in X$ and
$n>0$.
\end{prop}

\begin{proof}
 Since $\hat{\cz}(Q,\lambda)$ is a $U_\chi(\frak g)$-module,
  we know that its dual $\hat{\cz}(Q,\lambda)^{*}$ is a
  $U_{-\chi}(\frak g)$-module. Hence, we can get by \cite[Proposition 5.1]{FP1}
  $$\hbox{Ext}^n_{\mathcal{C}}(\hat{\cz}(Q,\lambda),
\hat{\cz}^{w^I}(Q,\mu^{w^I}))\simeq \hbox{H}^n(U_0(\frak g),
\hat{\cz}(Q,\lambda)^{*}\otimes\hat{\cz}^{w^I}(Q,\mu^{w^I}))$$
Noticing the fact $\Phi(||\frak g||_{\cz(Q,\lambda)})=
\Phi(||\frak g||_{\cz(Q,\lambda_i)^{*}})$, we can get that by
\cite[Proposition 6.2]{FP1}
$$\Phi(||\frak g||_{\cz(Q,\lambda)^{*}\otimes
\cz^{w^I}(Q,\mu^{w^I})})\subset \Phi(||\frak
g||_{\cz(Q,\lambda_i)})\cap \Phi(||\frak
g||_{\cz^{w^I}(Q,\mu^{w^I})})=0,$$
 the reason of
which is the same as in the proof of Theorem \ref{THMFORPRO}. Then
by Lemma \ref{FPTHM}. we can know that
$\hat{\cz}(Q,\lambda)^{*}\otimes\hat{\cz}^{w^I}(Q,\mu^{w^I})$ is a
projective $U_0(\frak g)$-module, as well as an injective module
because
 $U_0(\frak g)$ is a Frobenius algebra. Thus we have
 $\hbox{H}^n(U_0(\frak g),
\hat{\cz}(Q,\lambda)^{*}\otimes\hat{\cz}^{w^I}(Q,\mu^{w^I}))=0$
 and the lemma is proved.
 \end{proof}

\section{Andersen-Kaneda filtrations and sum formulas}

Maintain the notations as the previous sections. Especially, we
set $A=k[t]_{(t)}$  the localization of the polynomial ring $k[t]$
in one variable at the maximal ideal generated by $t$, we will
denote by $\tilde {A}$ the fraction field of $k[t]$.

  \subsection{Andersen-Kaneda filtrations}\label{subAKF} In the situation of \cite[\S 3]{AK}, Andersen and Kaneda
 constructed a filtration of the vector
space $F_{\lambda}(Q)=\hbox{Hom}_{G_1T}({{Z}}(\lambda)^\tau ,Q)$
for each projective $G_1T$-module $Q$, where $\tau$ denotes the
contra-variant dual \cite[\S 1.6]{AK}. We find that there exist
the similar filtration in the representation theory of modular Lie
algebra of $p$-character $\chi$ when $\chi$
 has standard Levi form. We will define this filtration in our situation and study its
 properties analogous to that in \cite[\S 3]{AK}.

Let $\hat Q$ be the projective module in category $\mathcal{C}_A$.
As stated in Remark \ref{PROJ}, there is a unique projective module
$\hat Q_k$ in category $\cc$ with
 $\hat Q_k\otimes_Ak\simeq \hat Q$.
Recall that $\hbox{Hom}_{\mathcal{C}_A}
(\hat{{Z}}_{A}^{w^I}(\lambda^{w^I}),\hat{{Z}} _A(\lambda))\cong A$
(Lemma \ref{LEMISO3.11}(2)). And there is a unique generator
$c:=\varpi'$ as appearing before Corollary \ref{KEYCOR}, up to units
in $A$, in $\hbox{Hom}_{\mathcal{C}_A}
(\hat{{Z}}_{A}^{w^I}(\lambda^{w^I}),\hat{{Z}} _A(\lambda)).$

Define
 \begin{eqnarray}\label{EQU4.1.1} & F^{\lambda}_A(\hat Q)
 =\hbox{Hom}_{\mathcal{C}_A}(\hat{{Z}}^{w^I}_A(\lambda^{w^I}),\hat
 Q),\nonumber\\
&E^{\lambda}_A(\hat Q)= \hbox{Hom}_{\mathcal{C}_A}(\hat
Q,\hat{{Z}}_A(\lambda)).
\end{eqnarray}
 By Lemmas \ref{ATILTINGDUAL}, Corollary \ref{PROJC} and Remark \ref{PROJ},  it's not hard to see that
\begin{equation} \label{EQU4.1.2}
F^{\lambda}_A(\hat Q)\mbox{ and } E^{\lambda}_A(\hat Q)\mbox{ are
both $A$-free module with the same rank, say } n_\lambda.%:=(\hat Q:\hat{{Z}}_A(\lambda))
\end{equation}
%where $(\hat Q_A:\hat{{Z}}_A(\lambda))$ denotes the multiplicity of
%$\hat{{Z}}_A(\lambda)$ in the $\hat{{Z}}$- filtration of $\hat Q.$
%In fact, ????????????????????//

%\begin{rem} Thus, from Proposition \ref{THMFORPRO} we know that each
%projective object in category $\mathcal C_A$ has a
%$\hat{\cz}_A$-filtration . Let $\hat{Q}_A(\lambda)$ be the
%projective cover of $\hat{L}_A(\lambda)$. Then the multiplicity of
%$\hat{{Z}}_A(\mu)$ in the $\hat{{Z}}$- filtration of $\hat
%Q_A(\lambda)$ is $$| W_I|\cdot
%[\hat{{Z}}_A(\mu):\hat{{L}}_A(\mu)]$$
% where $|W_I|$ is the number of elements in $W_I$ and
% $[\hat{{Z}}_A(\mu):\hat{{L}}_A(\lambda)]$ is the multiplicity of
% $\hat{{L}}_A(\lambda)$ as a composition factor of $\hat{{Z}}_A(\mu).$
%\end{rem}
Set $F^{\lambda}_k(\hat {Q}_k)=\hbox{Hom}_{\mathcal{C}}
(\hat{{Z}}^{w^I}_{\chi}(\lambda^{w^I}),\hat {Q}_k)$
 and $E^{\lambda}_k(\hat {Q}_k)=\hbox{Hom}_{\mathcal{C}_k}(\hat
{Q}_k,\hat{{Z}}_{\chi}(\lambda))$. The projectiveness of $\hat{Q}$
implies  by \cite[Proposition 3.3]{AJS} that
\begin{equation} \label{EQU4.1.3}
F^{\lambda}_k(\hat {Q}_k)\simeq F^{\lambda}_A(\hat Q)\otimes_Ak,
\end{equation}
and $E^{\lambda}_k(\hat {Q}_k)\simeq E^{\lambda}_A(\hat
Q)\otimes_Ak$.

 Similar to that in \cite[\S 3]{AK}, define a filtration
 $\{F^{\lambda}_A(\hat Q)^{(j)}\}_{j\geq 0}$ by setting
\begin{equation}\label{EQU4.1.4}F^{\lambda}_A(\hat Q)^{(j)}=\{\varphi\in F^{\lambda}_A(\hat
Q) \mid \psi\circ\varphi\in At^jc ,\psi\in E^{\lambda}_A(\hat Q)\}
\end{equation}
and let $F^{\lambda}_k(\hat {Q}_k)^{(j)}$ be the image of
$F^{\lambda}_A(\hat Q)^{(j)}$ in $F^{\lambda}_k(\hat {Q}_k)$, then
\begin{equation}\label{EQU4.1.5}
F^{\lambda}_k(\hat {Q}_k)^{(j)}\simeq(F^{\lambda}_A(\hat
Q)^{(j)}+tF^{\lambda}_A(\hat Q)) /tF^{\lambda}_A(\hat Q)\simeq
F^{\lambda}_A(\hat Q)^{(j)}/tF^{\lambda}_A(\hat Q)^{(j-1)}.
\end{equation}

\vskip0.2cm

We call such a filtration of the projective module $Q$ an {\sl{
Andersen-Kaneda filtration (or  AK filtration).}}

%Finally, we let $F^{\lambda}_k(\hat Q)_j$ denote the j-th quotient
%in the filtration $F^{\lambda}_k(\hat Q)=F^{\lambda}_k(\hat
%Q)^0\supset F^{\lambda}_k(\hat Q)^1\supset \cdots\supset
%F^{\lambda}_k(\hat Q)^r=0$, having
%\begin{equation}\label{EQU4.1.6}
%F^{\lambda}_k(\hat Q)_j=F^{\lambda}_k(\hat
%Q)^{(j)}/F^{\lambda}_k(\hat Q)^{j+1} \simeq F^{\lambda}_A(\hat
%Q)^{(j)}/(F^{\lambda}_A(\hat Q)^{j+1}+tF^{\lambda}_A(\hat
%Q)^{(j-1)}),j=0,1, \cdots,r-1\end{equation}

%\begin{rem} Similarly, we can define the corresponding filtration
%of $\widetilde{F^{\lambda}_A(\hat Q)}=
%\hbox{Hom}_{\mathcal{C}_A}(\hat{{Z}}_A(\lambda) ,\hat Q_A)$.
%\end{rem}
\vskip0.2cm

 Thanks to Lemma \ref{HOMZZ}(2),
we have the pairing $$a_{\lambda}: F_A^{\lambda}(\hat Q) \times
E_A^{\lambda}(\hat Q)\longrightarrow A$$ given by
$\psi\circ\varphi=a_{\lambda}(\varphi,\psi)c,\psi\in
E_A^{\lambda}(\hat{Q}),\varphi\in F^{\lambda}_A(\hat{Q}).$ When
tensored with $\tilde{A}$ we  have a bilinear $\tilde A$-form
$(-,-)$ arising from this pairing. This $\tilde A$-bilinear form
is by definition non-degenerate. Furthermore, there is with the
pairing, an $A$-homomorphism
$$\theta_{\lambda}:F^{\lambda}_A(\hat{Q})\longrightarrow E^{\lambda}
_A(\hat Q)^{\vee}=\hbox{Hom}_A( E^{\lambda}_A(\hat Q),A)$$ defined
by
$$\theta_{\lambda}(\varphi):\psi\mapsto a_{\lambda}(\varphi,\psi). $$
From the non-degeneracy of the $\tilde A$-bilinear form,
$\theta_{\lambda}$ is an $\tilde A$-isomorphism, when tensorred
with $\tilde A$. There are some basic facts with $\theta$ as
follows.

\begin{lem}\label{DUALS} \label{subsection4.2} There exist bases in $F_A^{\lambda}(\hat Q)$ and $E_A^{\lambda}(\hat
Q)$: $\{f_1,f_2,\cdots f_{n_{\lambda}}\} ,\{e_1,e_2,\cdots
e_{n_{\lambda}}\}$ respectively, together with a sequence of
positive integers  $\{
m_{\lambda}(1),m_{\lambda}(2),\cdots,m_{\lambda}(n_{\lambda})\}$
such that
$$\theta_{\lambda}(f_i)=t^{m_{\lambda}(i)}e_i,i=1,2,\cdots,n_{\lambda}\mbox{.}$$
Moreover,
$$\sum_{j\geq1}\dim F_k^{\lambda}(\hat {Q}_k)^{(j)}=\nu_t\det(\theta_{\lambda})=\sum^{n_{\lambda}}_{i=1}
m_{\lambda}(i)$$
\end{lem}

\begin{proof} Note that both $F^\lambda_A(\hat Q)$ and $E^\lambda_A(\hat Q)$
are $A$-free of rank $n_\lambda$, and $\theta_\lambda$ is an $\tilde
A$-isomorphism. For a given basis $\{e_1,\cdots,e_{n_\lambda}\}$ in
$E^\lambda_A(\hat Q)$, there  must exist a basis
$\{f_1,\cdots,f_{n_\lambda}\}$ in $F^\lambda_A(\hat Q)$ such that
$$\theta_{\lambda}(f_i)=t^{m_{\lambda}(i)}e_i,i=1,2,\cdots,n_{\lambda}\mbox{.}$$
As to the second formula, it follows from the standard arguments
as in \cite[II \S 8.18]{Jan3}.
\end{proof}

%
%4.4.
\subsection{Connection with Jantzen  filtrations}\label{4.4} Recall there are subspace filtrations on
Weyl modules in representations of algebraic groups (cf.
\cite{Jan3}). It's usually called Jatzen's filtration. The idea can
be adopted to the Lie algebra case when $\chi$ is of standard Levi
forms. Recall there are up to units in $A$, unique generators
$c:=\varpi'$ and $c':=\varpi$ respectively in
$\Hom_A(\hat{{Z}}_A^{w^I}(\lambda^{w^I}), \hat{Z}_A(\lambda))$  and
in $\Hom_A(\hat{Z}_A(\lambda), \hat{{Z}}_A^{w^I}(\lambda^{w^I}))$
(see the first subsection). According to Corollary \ref{KEYCOR}, we
have
 \begin{equation}\label{EQU4.4.9}c\circ c'=t^{N(I,\lambda)}
\hbox{id}_{\hat{{Z}}_A^{w^I}(\lambda^{w^I})}.
\end{equation}

In our case, the Jantzen filtration on
$\hat{{Z}}_A^{w^I}(\lambda^{w^I})$ and $Z_A(\lambda)$ can be defined
some sequences of their vector subspaces respectively defined via
$$\aligned
&\hat{{Z}}_A^{w^I}(\lambda^{w^I})^{(j)}=\{v\in\hat{{Z}}_A^{w^I}
(\lambda^{w^I})\mid cv\in t^j\hat{{Z}}_A(\lambda)\}\cr
&\hat{{Z}}_A(\lambda)^{(j)}=\{v'\in\hat{{Z}}_A(\lambda)\mid c'v'\in
t^{j}\hat{{Z}}_A^{w^I}(\lambda^{w^I})\}.\endaligned $$
 (cf. \cite[3.8]{Jan2}).

 As argument in the proof of Lemma \ref{FORDUALB} and Remark \ref{FORDUALB2},
 it's not hard to see there exist bases
$\{v_1,v_2,\cdots,v_n\}$ of $\hat{{Z}}_A^{w^I}(\lambda^{w^I})$
 and bases $\{v_1',v_2',\cdots,v_n'\}$ of  $\hat{{Z}}_A(\lambda)$
  and integers
$a_1,a_2,\cdots,a_n\in\mathbb{N}$ such that

\begin{equation}\label{EQU4.4.10}cv_i=t^{a_i}v_i',\qquad\qquad\qquad
c'v_i'=t^{N(I,\lambda)-a_i}v_i.
\end{equation}

We denote by $\hat{{Z}}_{\chi}^{w^I}(\lambda^{w^I})^{(j)}
 (\mbox{resp. }\hat{{Z}}_{\chi}(\lambda))$ the image of
$\hat{{Z}}_A^{w^I}(\lambda^{w^I})^{(j)}$ in
 ${Z}_{\chi}^{w^I}(\lambda^{w^I})$
(resp. the image of $\hat{{Z}}_A(\lambda)^{(j)}$ in
$\hat{{Z}}_{\chi}(\lambda)$ ), then we have (set
$\bar{v_i}=v_i\otimes1,\bar{v_i}'=v_i'\otimes1 $)

\begin{eqnarray}\label{EQU4.4.11}
&\hat{{Z}}_{\chi}^{w^I}(\lambda^{w^I})^{(j)}=\sum_{i\atop a_i\geq
j}k\bar{v}_i, \nonumber\\
&\hat{{Z}}_{\chi}(\lambda)^{(j)}=\sum_{i\atop
{N(I,\lambda)-a_i}\geq j}k\bar{v}_j'
\end{eqnarray}

 Observing that
$\bar{v_i}$ and $\bar{v}_i'$ have the same weight we deduce from
(\ref{EQU4.4.11})

\begin{equation}\label{EQU4.4.12}
\hbox{ch}\hat{{Z}}_{\chi}^{w^I}(\lambda^{w^I})^{(j)}+
\hbox{ch}\hat{{Z}}_{\chi}(\lambda)^{(N(I,\lambda)-j+1)}
=\hbox{ch}\hat{{Z}}_{\chi}(\lambda)
\end{equation}

%4.5
%\subsection{Further description of AK filtrations}\label{4.5}
Let's return to the AK filtration.  With the above arguments, we
know $\hat Q_A$ has a $\hat{\cz}_A^{w^I}$-filtration,
$$\hat Q_A=\hat Q_A^{(0)}\supset \hat Q_A^{(1)}\supset\cdots\supset \hat Q_A^{(r)}=0$$
with

\begin{equation}\label{EQU4.2.7}\hat Q_A^{(i-1)}/\hat Q_A^{(i)}\simeq
\hat{{Z}}^{w^I}_A(\lambda_i^{w^I})^{n_i}, n_i=(\hat
Q_A:\hat{{Z}}_{A}(\lambda_i))
\end{equation}
satisfying as in Remark \ref{rem3.15} that
 \begin{equation}\label{EQU4.2.8}\mbox{if}\mbox{ }
\lambda_i+\mathbb{Z}I\leq\lambda_j+\mathbb{Z}I\mbox{, then}\mbox{
}i\geq j
\end{equation}
for $i=1,2,\cdots,r$.

 We can get an exact
sequence
\begin{eqnarray}\label{STAR}
0\rightarrow \Hom_{\mathcal{C}_A}
(\hat{{Z}}^{w^I}_A(\lambda_i^{w^I}),\hat Q_A^{(i-1)}) &\rightarrow
\Hom_{\mathcal{C}_A}(\hat{{Z}}_A^{w^I}(\lambda_i^{w^I}),
\hat{{Z}}^{w^I}_A(\lambda_i^{w^I})^{n_i}) \nonumber\\
 &\rightarrow
\hbox{Ext}^1_{\mathcal C_A}(\hat{{Z}}_A^{w^I}(\lambda_i^{w^I}),
\hat{Q}_A^{(i)})\rightarrow 0.
\end{eqnarray}
which follows from the exact sequence:
$$0\rightarrow \hat
Q_A^{(i)}\rightarrow \hat{Q}_A^{(i-1)} \rightarrow
\hat{{Z}}_A^{w^I}(\lambda_i^{w^I})^{n_i} \rightarrow0$$ operated by
the functor
$\hbox{Hom}_{\mathcal{C}_A}(\hat{{Z}}_A^{w^I}(\lambda_i^{w^I}),-)$,
% Then, by Remark \ref{rem3.17}
as well as the fact that
$$\hbox{Hom}_{\mathcal{C}_A}(\hat{{Z}}_A^{w^I}(\lambda_i^{w^I}),\hat
Q_A^{(i)})=0.$$
 Furthermore, the inclusion $\hat
Q_A^{(i-1)}\hookrightarrow \hat Q_A$ induces an isomorphism
\begin{eqnarray}\label{SIMEQ}\Hom_{\mathcal{C}_A}(\hat{{Z}}_A^{w^I}(\lambda_i^{w^I}),\hat
Q_A^{(i-1)})\simeq
\Hom_{\mathcal{C}_A}(\hat{{Z}}_A^{w^I}(\lambda_i^{w^I}), \hat
Q_A)=F_A^{\lambda_i}(\hat Q), \end{eqnarray} because $\hat
Q_A/\hat Q_A^{(i-1)}$ is filtered by
$\hat{{Z}}_A^{w^I}(\lambda_j^{w^I}),j>i$ .

Set $\widetilde{F^{\lambda}_A(\hat
Q)}=\hbox{Hom}_{\mathcal{C}_A}(\hat{{Z}} _A(\lambda),\hat Q_A)$.
Then it's not hard to see the following fact by summing up the
above arguments (the detailed proof may be referred to
\cite[3.5]{AK}).

\begin{lem}\label{lem4.5}
There exist $A$-bases $\{\psi_1,\psi_2,\cdots,\psi_r\} \mbox{ of }
E_A^{\lambda}(\hat Q)$ and $\{\psi'_1,\psi'_2,\cdots,\psi'_r\},$
of
 $\widetilde{F^{\lambda}_A(\hat Q)}
$ such that
$$\psi_i\circ\psi'_j=\delta_{ij}t^{N(I,\lambda)}
\hbox{id}_{\hat{{Z}}_A(\lambda)}, i,j=1,2,\cdots,r.\qed$$
\end{lem}

%\begin{proof} With the notation as in \S \ref{subAKF} we see from Lemma \ref{HOMZZ}
%that if $\lambda=\lambda_i$, then we have an isomorphism
%$$\widetilde{F^{\lambda}_A(Q)}\tilde{\longleftarrow}Hom_{\mathcal{C}_A}(\hat{{Z}}_A(\lambda),
%\hat Q_A^{(i-1)})\tilde{\longrightarrow}
%\hbox{Hom}_{\mathcal{C}_A}(\hat{{Z}}_A(\lambda),
%\hat{{Z}}_A^{w^I}(\lambda^{w^I})^{n_i})$$ induced by the inclusion
%$\hat Q_A^{(i-1)}\hookrightarrow \hat Q_A$ and projective $\hat
%Q_A^{i-1}\rightarrow \hat Q_A^{(i-1)}/\hat Q_A^i\cong
%\hat{Z}_A^{w^I}(\lambda^{w^I})^{n_i}$ respectively. Now let
%$\psi'\in \widetilde{F^{\lambda}_A(\hat Q)}$ be the element which
%under these isomorphisms corresponds to the composite of $c'$ with
%the j-th inclusion $i_j:
%\hat{{Z}}_A^{w^I}(\lambda^{w^I})\rightarrow
%\hat{{Z}}_A^{w^I}(\lambda^{w^I})^{n_i}$.
% Similar we have isomorphisms $$E^{\lambda}_A(\hat
%Q)\tilde{\longrightarrow} \hbox{Hom}_{\mathcal{C}_A}(\hat
%Q_A^{(i-1)},\hat{{Z}}_A(\lambda))\tilde{\longleftarrow}
%\hbox{Hom}_{\mathcal{C}_A}(\hat{{Z}}_A^{w^I}(\lambda^{w^I})
%^{n_i},\hat{{Z}}_A(\lambda))
%$$
%and we let $\psi_j\in E^{\lambda}_A(\hat Q)$ denote
% the element corresponds to the composite of the j-th projection
%$\pi_j:\hat{{Z}}_A^{w^I}(\lambda^{w^I})^{ n_i}\rightarrow
%\hat{{Z}}_A^{w^I}(\lambda^{w^I})$ with $c$.

%We get
%$$\psi_s\circ\psi_t'=(c\circ\pi_s)\cdot(i_t\circ c')
%=\delta_{st}c\circ c',s,t=1,2,\cdots,n_i$$ The lemma therefore
%follows from (\ref{EQU4.4.9}).
%\end{proof}
Recall for $\varphi\in F^{\lambda}_A(\hat Q)$ and $\psi\in
E^{\lambda}_A(\hat Q)$, we have
 $\psi\circ\varphi=a_{\lambda}(\varphi,\psi)c$. Hence $c\circ c'=t^{N(I,\lambda)}\id_{\hat{Z}^{N+1}_A}$ implies
implies
\begin{equation}\label{EQU4.6.13}
\mbox{}\psi\circ\varphi\circ c'=
a_{\lambda}(\varphi,\psi)t^{N(I,\lambda)}\id_{\hat{{Z}}_A(\lambda)}
\end{equation}
With the bases from the above lemma, we may write $\varphi\circ
c'=\sum^r_{s=1}b_s\psi'_s$ for some $b_s\in A$. Furthermore,
$b_s=a_{\lambda}(\varphi,\psi_s)$ for all $s$. Hence $\varphi\circ
c'\in t^j\widetilde{F^{\lambda}_A(\hat Q)}$ if and only if
$t^j\mid a_{\lambda}(\varphi,\psi_s)$ for all $s$, i.e., if and
only if $\varphi \in F^{\lambda}_A(\hat Q)^{(j)}$.

Thus,  we have  a corollary to the above lemma.
\begin{lem}\label{prop4.6}
For each $j\in\mathbb{N}$ we have
$$F_A^{\lambda}(\hat Q)^{(j)}=\{\psi\in F^{\lambda}_A(\hat Q)\mid\varphi\circ
c' \in t^j\widetilde{F^{\lambda}_A(\hat Q)}\}.\qed$$
\end{lem}

 %4.7.
\subsection{Sum formulas over $\cc$}
Now we are in the position to present the main results about AK
filtrations and the related sum formulas.

\begin{lem} \label{LASTLEM}
%\begin{itemize}
%\item[(1)] Let $\hat{Q}_A$ be a module in ${\cc_A}$ and
%$\hat{Q}=\hat{Q}_A\otimes_Ak$. Then $\Hom_{\cc_A}(\hat{Z}_A^i,\hat
%Q)\otimes_A k\cong \Hom_\cc(\hat{Z}_\chi^i,Q)$.
% \item[(2)]
$\Ker\bar\varphi_i\cong \Coker\bar\varphi_i$, where $\bar
 \varphi_i$ is an extension of $\varphi_i$ by change of rings
 $-\otimes_Ak$,
 and $\varphi_i$ is defined as in Lemma \ref{AK1.9}.
%\end{itemize}
\end{lem}

\begin{proof} %(1) Recall each object in $\cc_A$ is required to be finitely generated
%as $A$-module. So the statement follows from \cite[Theorem
%2.38]{CR}.
It follows from  \cite[Lemma 3.5]{Jan2}.
\end{proof}

 \begin{thm}\label{SumFor1}
Let $\nu\in X$. Then
$$\sum_{j\geq1}\hbox{\dim}F^{\lambda}_k(\hat Q_\chi(\nu))^{(j)}
=\sum_{j\geq1}[\hat{{Z}}_\chi(\lambda)^{(j)}:\hat{{L}}_\chi(\nu)]$$
\end{thm}

\begin{proof}

Let $\hat Q=\hat Q_\chi(\nu)\in\mathcal{C},\hat Q_A=\hat Q_A(\nu)\in
\mathcal {C_A}$. Denote by $\phi_i:
\Hom_{\mathcal{C}_A}(\hat{{Z}}_A^{i+1},\hat Q_A)\rightarrow
\Hom_{\mathcal{C}_A}(\hat{{Z}}_A^{i},\hat Q_A)$ which is induced by
$\varphi_i$. Here $\varphi_i$ is  defined as in Lemma \ref{AK1.9}.
Let $\phi=\phi_1\circ\cdots\circ\phi_N$. Then Lemma \ref{prop4.6}
says
$$F^{\lambda}_A(\hat Q)^{(j)}=\{\varphi\in
F^{\lambda}_A(Q)\mid \phi(\varphi)\in
t^j\widetilde{F^{\lambda}_A(Q)}\}.$$
Analogy of the arguments as
in \S\ref{subsection4.2} gives the following sum formula.
 \begin{equation}\label{EQU4.5.14}\sum_{j\geq1}\hbox{dim}F^{\lambda}_k(\hat
Q)^{(j)}=\nu_t(\hbox{det}\phi)
=\sum^N_{i=1}\nu_t(\hbox{det}\phi_i).
\end{equation}

Set $C_i=\Coker\varphi_i$. Then by \cite[\S 3.7]{Jan2}, we have

\begin{equation}\label{EQU4.5.15}tC_i=0\mbox{, for all i}.
\end{equation}

 Hence we also have

\begin{equation}\label{EQU4.5.16}
t\hbox{Ext}_{\mathcal{C}_A}^1(C_i,\hat Q_A)=0\mbox{, for all i}.
\end{equation}

Notice that  the homomorphism
$\hat{{Z}}_A^{i-1}\rightarrow\hat{{Z}}_A^i$ is injective,
%\cite[Remark 3.6]{Jan2},
we have the short exact sequence

$$0\rightarrow\hat{{Z}}_A^{i}\stackrel{\varphi_{i}}{\rightarrow}
\hat{{Z}}_A^{i+1}\stackrel{\tilde{\varphi_i}}{\rightarrow}C_i\rightarrow0.
$$

Moreover we can get the exact sequence through the action by
$\Hom_{\mathcal{C_A}}(-,\hat Q_A)$

\begin{equation}\label{EQU4.5.17}
0\rightarrow \hbox{Hom}_{\mathcal{C_A}}(\hat{{Z}}_A^{i+1},\hat
Q_A)\stackrel{\phi_i}{\rightarrow}
\hbox{Hom}_{\mathcal{C_A}}(\hat{{Z}}_A^{i},\hat Q_A)\rightarrow
 \hbox{Ker}\Phi_i\rightarrow0
 \end{equation}
where $\Phi_i:\Ext_{\mathcal{C}_A}^1(C_i,\hat Q_A)\rightarrow
\hbox{Ext}_{\mathcal{C}_A}^1(\hat{{Z}}_A^{i+1},\hat Q_A)$. It is
clear that $\Ker\Phi_i$ is the submodule of
$\hbox{Ext}_{\mathcal{C}_A}^1(C_i,\hat Q_A)$. Hence
(\ref{EQU4.5.15}) implies

\begin{equation}\label{EQU4.5.18}
\nu_t(\hbox{det}\phi_i)
%=l_t(\Coker\phi_i)
=\hbox{dim}_k(\hbox{Ker}\Phi_i\otimes_Ak)
\end{equation}

Set $\bar{\varphi_i}=\varphi_i\otimes_Ak,\bar{C_i}=C_i\otimes_Ak$.
Then $\bar{C_i}=\Coker\bar{\varphi_i}$ and the injectivity of
$\hat Q$ gives the exactness of the top row in the following
diagram (note that $U_\chi(\ggg)$ is a Frobenius algebra,  the
projective $U_{\chi}(\frak g)$-modules is injective): \vskip5pt
$${\tiny
\begin{diagram} 0 & \rTo & \hbox{Hom}_\cc(\bar C_i,\hat Q)& \rTo &
\hbox{Hom}_\cc(\hat{{Z}}_\cc^{i+1},Q) & \rTo &
\hbox{Hom}_\cc(\hat{{Z}}_\cc^{i},\hat Q)&\rTo&
\hbox{Hom}_\cc(\hbox{Ker}
\varphi_i,Q)&\rTo&0\\
&&&&\dTo_{\simeq}&&\dTo_{\simeq}&&&&\\
&&&&\hbox{Hom}_{\mathcal{C}_A} (\hat{{Z}}_A^{i+1}, \hat
Q_A)\otimes_Ak &\rTo& \hbox{Hom}_{\mathcal{C}_A}(\hat{{Z}}_A^{i},
\hat Q_A)\otimes_Ak&\rTo& \Ker\Phi_i\otimes_Ak&\rTo&0
\end{diagram}}$$

\noindent Here the bottom sequence is obtained by tensoring
(\ref{EQU4.5.17}) with $k$ and the vertical isomorphism comes from
the application of Proposition \ref{THMFORPRO} and Lemma \ref{HOMZZ}
along with Remarks \ref{PROJ}(2), \ref{REMDUAL} and \ref{PROJCORA}.
The diagram shows that
 $\hbox{Ker}\Phi_i\otimes_Ak$ may be identified with
$\hbox{Hom}_{\mathcal{C}}(\hbox{Ker}\bar{\varphi_i},\hat Q)$,
thereby  this has the same dimension as
$\hbox{Hom}_{\mathcal{C}}(\bar{C_i},\hat
 Q)$, owing to Lemma \ref{LASTLEM}.
  By  Lemma \ref{TILTINGDUAL}
  and Remark \ref{REMDUAL},  it can be seen that this dimension equals the multiplicity
  $[\bar C_i, \hat{L}_\chi(\mu)]$.
  Consider $\bar C_i=\Ker \bar \varphi_i$ (Lemma \ref{LASTLEM}), which is decided by the
  formula in Remark \ref{FORDUALB2}. Hence, the arguments on (\ref{EQU4.4.10})-(\ref{EQU4.4.12}) tells us
  \begin{equation}\label{EQU4.5.19}
\sum_{j\geq1}\hbox{ch}\hat{{Z}}_{\chi}(\lambda)^{(j)}=
\sum^N_{i=1}\hbox{ch}\bar{C}_i.
\end{equation}
  Combining (\ref{EQU4.5.14}), (\ref{EQU4.5.18}) and
(\ref{EQU4.5.19}), we know
$$ \sum_{j\geq1}\hbox{dim}F^{\lambda}_k(\hat
Q)^{(j)}=\sum^N_{i=1} [\bar C_i: \hat{L}_\chi(\nu)] =\sum_{j\geq1}
[\hat{{Z}}_{\chi}(\lambda)^{(j)}:\hat{{L}}_{\chi}(\nu)].$$ The proof
is completed.
\end{proof}

%4.9
The following result shows the connection between the dimension of
individual terms between
 AK filtrations and Jantzen's filtrations.

\begin{thm}\label{SumFor2} Let $\hat Q\in\mathcal{C}$ be
projective. Then we have a formula $\dim F^{\lambda}_k(\hat Q)^{(j)}
=\dim \Hom(\hat{{Z}}_{\chi}(\lambda)^{(j)},\hat Q),\mbox{ for all
}j\in\mathbb{N}.$ Especially,  for $\hat Q=\hat Q(\nu),\nu\in X$,
$\label{EQU4.9.22} \dim F^{\lambda}_k(\hat Q(\nu))^{(j)}
=[\hat{{Z}}_{\chi}(\lambda)^{(j)}:\hat{{L}}_{\chi}(\nu)], \mbox{ for
all }j\mbox{. }$
\end{thm}

\begin{proof} The second part is a direct implication of the first
one. It's sufficient to prove the first one. For this, let us
verify the statement below.
\begin{equation}\label{EQU4.8.20}\mbox{If
}\bar{\varphi}\in F^{\lambda}_k(\hat Q)^{(j)},\mbox{ then }
\bar{\varphi}(\hat{{Z}}_{\chi}^{w^I}(\lambda^{w^I})^{(N(I,\lambda)+1-j)})=0.
\end{equation}

 We denote $\bar{\varphi}\in F^{\lambda}_k(\hat
Q)^{(j)}$ is the image of some $\varphi\in F^{\lambda}_A(\hat
Q)^{(j)}$ in
 $F^{\lambda}_k(\hat Q)$. By Proposition \ref{prop4.6} we have $\varphi\circ c'\in t^j\widetilde{
F^{\lambda}_A(\hat Q)}.$ On the other hand, by (\ref{EQU4.4.10})
we have
$$\varphi\circ
c'(v_i')=t^{N(I,\lambda)-a_i}\varphi(v_i),i=1,2,\cdots n$$
 We conclude
that
\begin{equation}\label{EQU4.8.21}\mbox{if }N(I,\lambda)-a_i\leq j\mbox{ then }\varphi(v_i)\in
t^{j-N(I,\lambda)+a_i} \hat Q_A
\end{equation}
 By (\ref{EQU4.4.11}) we
see that $\hat{{Z}}_A^{w^I}(\lambda^{w^I})^{(N(I,\lambda)+1-j)}$
is spanned by $\bar{v_i}$ where $a_i\geq N(I,\lambda)+1-j.$ Hence
(\ref{EQU4.8.21}) implies that if
$\bar{v_i}\in{Z}_A^{w^I}(\lambda^{w^I})^{(N(I,\lambda)+1-j)}$ ,
then $\varphi(v_i)\in t\hat Q_A.$ That is to say that
$\bar{\varphi}(\bar{v_i})=0$, thereby  (\ref{EQU4.8.20}) is
proved.

According to (\ref{EQU4.8.20}), we see that
$$F^{\lambda}_k(\hat Q)^j\subseteq
\Hom_{\mathcal{C}_A}(\hat{{Z}}_{\chi}^{w^I}(\lambda^{w^I})
/\hat{{Z}}_{\chi}^{w^I} (\lambda^{w^I})^{(N(I,\lambda)+1-j)},\hat
Q).$$ Since $\hat Q$ is injective, the dimension of this
homomorphism space only depends on the character of $
\hat{{Z}}_{\chi}^{w^I}(\lambda^{w^I})/
\hat{{Z}}_{\chi}^{w^I}(\lambda^{w^I})^{(N(I,\lambda)+1-j)}$ .
Therefore by (\ref{EQU4.4.12}) we deduce
$$\hbox{dim}F^{\lambda}_k(\hat Q)^{(j)} \leq \dim
\Hom_{\mathcal{C}_A}(\hat{{Z}}_{\chi}(\lambda),\hat Q)^{(j)} \mbox{
for all j}\in\mathbb{N}.$$ On the other hand, each
finite-dimensional projective object in $\cc$ is the  direct sum of
some $Q(\nu)$. Thus Theorem \ref{SumFor1} make in force the equality
true. The proof is completed.
\end{proof}

\begin{rem} %If  $\hat Q\in\mathcal{C}$ is projective, we denote
%for $\nu\in X$ by $(\hat Q:\hat Q(\nu))$
 %the number of times $\hat Q(\nu)$ occurs as summand in $\hat Q$. We set $R^+(\lambda)=
 %\{\lambda\in R^+\mid\langle\lambda+\rho,\alpha^{\vee}\rangle\not\equiv 0\;(\hbox{mod }p)\}$
%and we let $n_{\alpha}\in\{1,\cdots,p-1\},\alpha\in R^+(\lambda)$
%denote the residue
% $\langle\lambda+\rho,\alpha^{\vee}\rangle \;(\hbox{mod }p)$.
For each projective module $\hat Q\in\mathcal{C}$ we have
\begin{itemize}
 \item[(1)] The length of the  Jantzen
filtration of $\hat{{Z}}_{\chi}(\lambda)$ is just $N(I,\lambda)$.
So it follows from the above theorem that
$$ F^{\lambda}_k(\hat
Q)^{(N(I,\lambda)+1)}=0 \mbox{ for all projective modules }\hat
Q\in\mathcal{C}\mbox{.}
$$
 \item[(2)]  We have another version of the sum
formula
\begin{eqnarray*}\sum_{j\geq1}\hbox{dim}F^{\lambda}_k(\hat
Q^{(j)})&=&\sum_{\alpha\in R^+(\lambda)\backslash R^+_I}
(\sum_{i\geq0}(\hat Q:\hat{{Z}}_A(\lambda-(ip+n_{\alpha})))-\\
&&\sum_{i\geq1}(\hat Q: \hat{{Z}}_A(\lambda-ip\alpha))).
\end{eqnarray*}
\item[(3)]  $\dim F^{\lambda}_k(\hat Q)_0=(\hat Q:\hat
Q(\lambda))$ and
 $\dim F^{\lambda}_k(\hat Q)_{N(I,\lambda)}=(\hat
Q:\hat Q(\lambda'))\mbox{ where }
\hat{{L}}(\lambda')=\hbox{Soc}(\hat{{Z}}_{\chi}(\lambda)).$
\end{itemize}
\end{rem}


\begin{thebibliography}{9}

\bibitem{AJS}Andersen H. H., Jantzen J. C. and Soergel, W., {\em Representations of quantum groups at a pth
 root of unity and of semisimple groups in charactertic p: independence of p.}, Asterisque 220 (1994).

\bibitem{AK} Andersen H. H. and Kaneda, M.,  {\em Filtrations on $G_1$T-modules},
Proc. Londen Math. Soc. 82 (2001), 614-646.

%\bibitem{An98}H. H. Andersen, {\em A sum formular for tilting filtration}, Proceedings of a conference
% in honor of D. Buchsbaum, Rome 1998, J. Pure Appl. Algebra 152 (2000), 17-40.

\bibitem{CPS} Cline E., Parshall B. and Scott, L.,  {\em On injective modules for
 infinitesimal algebraic groups I}, J. London Math. Soc.(2) 31 (1985), 277-291

\bibitem{CR} Curtis C.W. and Reiner I., {\em Methods of
representation theory. With applications to finite groups and
orders}, Vol.I, A Wiley-Interscience Publication, New York, 1981.

\bibitem{FP1} Friedlander E. M. and Parshall B., {\em Modular representation theory of Lie
 algebras}, Amer. J. Math. 110 (1988), 1055-1093.

\bibitem{FP2} Friedlander E. M. and Parshall B., {\em Geometry of p-unipotent Lie algebras}, J. Algebra
 109 (1987), 25-45

\bibitem{FP3} Friedlander E. M. and Parshall B., {\em Support varieties for restricted Lie algebra},
Invent. math. 86 (1986),553-562

\bibitem{FP4} Friedlander E. M. and Parshall B., {\em Deformations
of Lie algebra representations}, Amer. J. Math. 112(1990), 375-395.

%\bibitem{Hum1} J. Humphreys, {\em Introduction to Lie algebras and
%their representations},

%\bibitem{Hum2} J. Humphreys, {\em Linear algeberaic groups}

\bibitem{Hum3} Humphreys J.. {\em Modular representations of classical
Lie algebras and semisimple groups}, J. Algebra 19 (1970), 51-79


\bibitem{Jan1} Jantzen J. C., {\em Representations of Lie algebras in prime characteristic}, in ; A.
              Broer (Ed. ), Representation Theories and Algebraic Geometry, Proceedings,
              Montreal, 1997, in: NATO ASI Series, Vol. C514, Kluwer, Dordrecht, 1998, pp.
              185-235.

\bibitem{Jan2} Jantzen J. C.,  {\em Modular representations of reductive Lie algebras}, Journal of
               Pure and Applied Algebra 152 (2000), 133-185.

\bibitem{Jan3}  Jantzen J. C.,  {\em Reresentations of Algebraic Groups}, 2nd edn, American Mathematical
 Society, Providence, RI, 2003.

%\bibitem{Jan4} Jantzen J. C.,  {\em Subregular nilpotent representations of $\sss\lll_n$ and
  %$\sss\ooo_{2n+1}$}, Math. Proc. Cambridge Philos. Soc. 126 (1999), 223-257.

%\bibitem{Jan5} Jantzen J. C.,  {\em Representations of $\frak{s0}_5$ in prime characterstic}, Preprint
 % Aarhus series 1997:13.

\bibitem{Jan6} Jantzen J. C.,  {\em Kohomologie Von p-Lie-Algebra und nilpotente Elemente}, Abh. Math. Sem. Univ.
  Hamburg 56 (1986), 191-219.

\bibitem{Je} Jessen B. R., {\em Representation Theorey of Lie algebras in prime
characteristic}, Progress Report, Aarhus University, 1999.

\bibitem{L} Lusztig G., {\em Periodic $W$-graphs}, Represent.
Theory 1 (1997), 207-279.

\bibitem{KW}  Kac V. and Weisfeiler B., {\em Coadjoint action of a semisimple algebraic group and the center of
   the enveloping algebra in characteristic p}, Indag. Math. 38 (1976), 136-151.

%\bibitem{Pre}A. Premet, {\em Support varities of non-restricted modules over Lie algebras of reductive groups},
%J. London Math. Soc. (2)  55 (1997), 236-250

\end{thebibliography}
\end{document}